# A LARGE DEVIATIONS APPROACH TO ASYMPTOTICALLY OPTIMAL CONTROL OF CRISSCROSS NETWORK IN HEAVY TRAFFIC[1]

By Amarjit Budhiraja and Arka Prasanna Ghosh

*University of North Carolina at Chapel Hill*

In this work we study the problem of asymptotically optimal control of a well-known multi-class queuing network, referred to as the "crisscross network," in heavy traffic. We consider exponential inter-arrival and service times, linear holding cost and an infinite horizon discounted cost criterion. In a suitable parameter regime, this problem has been studied in detail by Martins, Shreve and Soner [*SIAM J. Control Optim.* **34** (1996) 2133–2171] using viscosity solution methods. In this work, using the pathwise solution of the Brownian control problem, we present an elementary and transparent treatment of the problem (with the *identical* parameter regime as in [*SIAM J. Control Optim.* **34** (1996) 2133–2171]) using large deviation ideas introduced in [*Ann. Appl. Probab.* **10** (2000) 75–103, *Ann. Appl. Probab.* **11** (2001) 608–649]. We obtain an asymptotically optimal scheduling policy which is of threshold type. The proof is of independent interest since it is one of the few results which gives the asymptotic optimality of a control policy for a network with a more than one-dimensional workload process.

## 1. Introduction.

Stochastic networks are ubiquitous in problems involving manufacturing, communication and computer systems. Designing good controls for general multi-class networks is an important and challenging problem. In recent years, using tools from diffusion approximations, there has been a significant progress in obtaining asymptotically optimal controls for a broad range of stochastic networks in heavy traffic. One common approach to the optimality question is via certain singular control problems,

Received April 2003; revised October 2004.

[1]Supported in part by ARO Grant W911NF-04-1-0230.

the so-called Brownian control problems (BCP), which are obtained as "formal" heavy traffic limits of queuing networks. There are several works (e.g., [6, 7]) which use the optimal solution of the BCP to construct control policies for the corresponding queuing networks. These policies seem to perform quite well in simulation studies, however, there are relatively few results showing asymptotic optimality of such policies. Recently, in [1, 4], using large deviation ideas, a promising technique for addressing asymptotic optimality questions has been introduced. Using these techniques, the authors prove asymptotic optimality of a certain threshold-based scheduling policy for a "parallel server model." Other recent results on asymptotic optimality of control policies for stochastic networks are in [8, 10, 11, 12, 13].

In the current work we study a well-known model, often referred to as the "crisscross network." It has been studied in [7, 16] and in great detail in [11] and [10]. The network is described in detail in Section 2. The basic problem is the optimal sequencing of jobs in a two station-two customer queuing system. We consider linear holding costs and an infinite horizon discounted cost criterion [see (2.17)]. We believe the scheduling policy that we propose will also be asymptotically optimal for a finite time horizon cost criterion with a linear holding cost. However, for the sake of simplicity, we restrict our attention to the first criterion. Even though the network is quite simple to describe, the analysis of the control problem is rather subtle in that the form of an asymptotically optimal scheduling policy and the methods of proof seem to strongly depend on the parameter regime under consideration. Broadly, one can divide the study of the problem into two different parameter regimes: *Case I*: $h_1\mu_1 - h_2\mu_2 + h_3\mu_2 \leq 0$ and *Case II*: $h_1\mu_1 - h_2\mu_2 + h_3\mu_2 > 0$, where $h_i$'s are the holding costs and $\mu_i$'s are the asymptotic service rates (see Section 2 for precise definitions). *Case I* yields a simple threshold policy and the proof of asymptotic optimality of this policy is given in [16].

*Case II* is the difficult case and in [11] its analysis has been subdivided into 4 subcases: *Case IIA*, ..., *Case IID*. In *Case IIA*, in addition to conditions of *Case II*, both $h_2\mu_2 - h_3\mu_2$ and $h_2\mu_2 - h_1\mu_1$ are nonnegative. The other three subcases (*Case IIB*, *Case IIC*, *Case IID*) correspond to either one of these two quantities being negative and the case where both are negative. Among the four subcases, *Case IIA* is most amenable to analysis, since in this case, the "effective cost" in the reduced workload formulation, $\hat{h}(w_1, w_2)$ is monotonic in both $w_1$ and $w_2$ (see Remark 3.3). This monotonicity is critical in obtaining an explicit, *pathwise* solution to the BCP. *Case IIA* was studied in [7] with specific numerical values of the parameters and though the authors did not prove asymptotic optimality of the proposed policy, they provided results from simulation studies indicating good performance of the control policy. This subcase was studied in complete detail in [11]. The authors proposed a control policy and proved the *near asymptotic* *optimality* of the policy by using quite technical machinery from viscosity solution analysis of Hamilton–Jacobi–Bellman (HJB) equations. The proposed policy is difficult to interpret and is not very intuitive. In addition, extension of the methodology to more complex networks appears quite daunting. One of the technical obstacles in such extensions is that a general theory of classical solutions to HJB equations or characterization results for the value function via viscosity solutions of HJB equations for the limiting control problem are not readily available. Finally, we note that, strictly speaking, [11] does not obtain an asymptotically optimal policy. By "near optimality," it is meant that for each $\eta > 0$, one can get a control policy (depending on $\eta$) which, asymptotically, is $\eta$-close to an asymptotically optimal strategy. In [10], using techniques from weak convergence theory, the authors show that the optimal costs for the queuing network problem can be well approximated by those for the optimization problem of the limiting control problem. The approach is quite general and powerful, but the authors do not obtain an actual control policy which is asymptotically (near) optimal.

In the current work we revisit the above problem (under the same parameter regime, namely, *Case IIA*) using a rather different approach introduced in the context of a "parallel server model" in [1, 4]. We present the BCP associated with this control problem and give the equivalent workload formulation. The BCP that we obtain is somewhat different from the one presented in [11]. Indeed, the authors there remark that their BCP is not well posed (see Remark 3.5 for more details on this). However, as we show in Section 3.1, the BCP presented in this paper has an explicit pathwise optimal solution. The scheduling policy we propose (see Definition 3.6) is directly motivated by the solution of the BCP, and therefore is easy to interpret. The scheduling policy is of threshold type and thus is quite simple to implement as well. In addition, our proof of the asymptotic optimality of the policy uses rather basic large deviation ideas which, we believe, can be extended to more general situations.

All inter-arrival and service times in this work will be assumed to be exponentially distributed. Proofs of many of the results in this paper can be extended to the case where the inter-arrival and service times are i.i.d. with distributions that satisfy suitable large deviation estimates. Indeed, in the parallel server model [1], the authors prove asymptotic optimality under precisely such assumptions on the underlying distributions. One important difference in our analysis is that in one of the key results of this paper (Theorem 4.9), in addition to the one-dimensional large deviation estimates that are crucially used in the proofs of [1], we also need *sample path* large deviation estimates (Theorem 5.1) for the underlying renewal processes. For a more detailed discussion on extending the results of this paper to the nonexponential case, see Remark 5.4.

The paper is organized as follows. The network is described in Section 2, along with the formulation of the problem and assumptions. In Section 3 we formulate the associated BCP and the corresponding equivalent workload formulation. We then propose a policy that is motivated by the equivalent workload formulation and the solution of the BCP. In Section 4 the asymptotic optimality of the proposed policy is proved through the two main results of the paper, Theorems 4.1 and 4.2. These results are as follows. Denoting the minimum cost associated with the BCP as $J^*$ and the cost associated with *any* control policy $T^r$ for the $r$th network as $\hat{J}^r(T^r)$, we show in Theorem 4.1 that

$$\liminf_{r \to \infty} \hat{J}^r(T^r) \geq J^*.$$

In Theorem 4.2 it is shown that, in the above display, the equality is achieved if $\{T^r\}$ is the sequence of policy proposed in Definition 3.6 of Section 3, with an appropriate choice of threshold parameters. The key steps in the proof of the two main theorems are in Theorems 4.8, 4.9 and 4.11. The proofs of these theorems are provided in Section 5.

## 2. The crisscross network.

2.1. *Queueing network model.* We consider a sequence of networks indexed by $r$, $r \in \mathbb{S} \subseteq \mathbb{R}^+$, where $\mathbb{S}$ is a countable set: $\{r_1, r_2, \ldots\}$ with $1 \leq r_1 < r_2 < \cdots$ and $r_n \to \infty$, as $n \to \infty$. A sketch of the $r$th network is described in Figure 1. Description for the $r$th network is as follows. For $i = 1, 2$, customers (or jobs) of Class $i$ arrive according to a Poisson process with rate $\lambda_i^r$ and have independent exponential service times at *Server* 1 with parameter $\mu_i^r$. Class 1 customers, once served by *Server* 1, leave the system. Class 2 customers, after being served by *Server* 1, proceed to *Buffer* 3 and are redesignated as Class 3 customers. There they are served by *Server* 2. They have i.i.d. exponential service times with parameter $\mu_3^r$. After service, these jobs exit the system. All inter-arrival and service times are assumed to be mutually independent and all buffers have infinite capacity. We also assume that the system starts empty.

2.2. *Preliminaries.* Let $(\Omega, \mathcal{F}, \mathbb{P})$ be a complete probability space. All the random variables and stochastic processes in this paper are assumed to be defined on this probability space. There is no loss of generality in making this assumption since we work with an expected loss function (see Section 2.4 for the definition of cost) and one can always enlarge the probability space to support all the processes considered in this paper. The expectation operation under $\mathbb{P}$ will be denoted by $\mathbb{E}$.

For each positive integer $m \geq 1$, let $\mathcal{D}^m$ be the space of right continuous paths with left limits, from $[0, \infty)$ to $\mathbb{R}^m$, with the usual Skorohod topology and let $\mathcal{B}(\mathcal{D}^m)$ be the corresponding Borel sigma-field. All of the continuous-time processes considered in this paper will have sample paths in $\mathcal{D}^m$. If $\{Z_n\}$ and $Z$ are processes with paths in $\mathcal{D}^m$ such that $Z_n$ converges weakly to $Z$ as $n \to \infty$, we will use the notation $Z_n \Rightarrow Z$ to denote this.

For each $r \in \mathbb{S}$ and $k = 1, 2$, let $\{u_k^r(i) : i = 1, 2, \ldots\}$ be a sequence of i.i.d. exponential random variables with mean $1/\lambda_k^r \in (0, \infty)$. We interpret $u_k^r(i)$ to be the time (in the $r$th network) between the arrival of the $(i-1)$st and the $i$th job for *Buffer* $k$ $(k = 1, 2)$. Similarly, the service times of the three different classes of jobs are defined as sequences of i.i.d. exponential variables $\{v_j^r(i) : i = 1, 2, \ldots\}$, with mean $1/\mu_j^r \in (0, \infty)$, $j = 1, 2, 3$, corresponding to the three classes. We also assume that the inter-arrival time sequence $\{u_k^r(i) : i = 1, 2, \ldots\}$, $k = 1, 2$, and the service time sequence $\{v_j^r(i) : i = 1, 2, \ldots\}$, $j = 1, 2, 3$, are mutually independent for each $r \in \mathbb{S}$.

Define

$$\xi_k^r(n) \doteq \sum_{i=1}^{n} u_k^r(i) \qquad \text{for } n = 1, 2, \ldots, k = 1, 2;$$

$$\eta_j^r(n) \doteq \sum_{i=1}^{n} v_j^r(i) \qquad \text{for } n = 1, 2, \ldots, j = 1, 2, 3.$$

The arrival and service processes are defined in terms of these as follows:

$$A_k^r(t) \doteq \sup\{n \geq 0 : \xi_k^r(n) \leq t\}, \qquad t \geq 0, k = 1, 2,$$

$$S_j^r(t) \doteq \sup\{n \geq 0 : \eta_j^r(n) \leq t\}, \qquad t \geq 0, j = 1, 2, 3.$$

The symbol $A_k^r(t)$ represents the number of jobs (customers) that have arrived in *Buffer* $k$ up to time $t$. The process $S_j^r(t)$ counts the number of jobs that *Server* $j$ could have completed if it had worked continuously during the interval $[0, t]$. Note that by our assumptions on the inter-arrival and service times, $A_k^r(\cdot)$ and $S_j^r(\cdot)$ are Poisson processes with rates $\lambda_k^r$ and $\mu_j^r$, respectively, for $k = 1, 2; j = 1, 2, 3; r \in \mathbb{S}$. For notational simplicity, throughout the paper, we will write the limit along the sequence $r_n$ as $n \to \infty$ simply as "$r \to \infty$." Also, $r$ will always be taken to be an element of $\mathbb{S}$ and, thus, hereafter, the qualifier $r \in \mathbb{S}$ will not be stated, explicitly. We assume that as $r \to \infty$, these rates approach finite limits, namely, we make the following assumption.

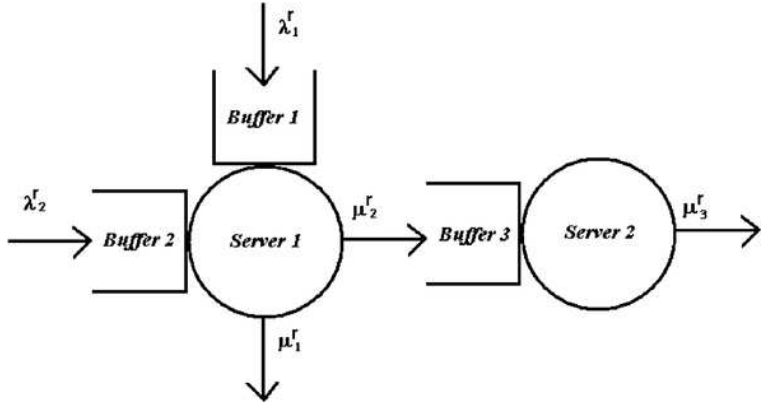

ASSUMPTION 2.1. There exist $\lambda_k \in (0, \infty), k = 1, 2$, and $\mu_j \in (0, \infty), j = 1, 2, 3$, such that

$$\lim_{r \to \infty} \lambda_k^r = \lambda_k, \qquad \lim_{r \to \infty} \mu_j^r = \mu_j, \qquad k = 1, 2; j = 1, 2, 3.$$

2.3. *Scheduling control.* Scheduling control for the $r$th network is described by a vector-valued service allocation process

$$T^r(t) \equiv (T_1^r(t), T_2^r(t), T_3^r(t)), \qquad t \geq 0,$$

where, for $j = 1, 2, 3$, $T_j^r(t)$ denotes the cumulative amount of service time devoted to activity $j$ (viz., working on Class $j$ jobs by the responsible server) in the time interval $[0, t]$. The idle-time processes are defined as follows:

$$I_1^r(t) \doteq t - T_1^r(t) - T_2^r(t), \qquad I_2^r(t) \doteq t - T_3^r(t), \qquad t \geq 0.$$

For $i = 1, 2$, $t \geq 0$, $I_i^r(t)$ represents the cumulative amount of time that the $i$th server has been idle in the time interval $[0, t]$. Recall that we assume that the system is initially empty. Thus, the three queue-length processes corresponding to the three buffers can be described as follows. For $t \geq 0$,

(2.1)
$$\begin{aligned} Q_i^r(t) &= A_i^r(t) - S_i^r(T_i^r(t)), \qquad i = 1, 2, \\ Q_3^r(t) &= S_2^r(T_2^r(t)) - S_3^r(T_3^r(t)). \end{aligned}$$

The workload process $W^r(\cdot) = \{(W_1^r(t), W_2^r(t)), t \geq 0\}$ is defined as follows. For $t \geq 0$,

(2.2)
$$\begin{aligned} W_1^r(t) &\doteq \frac{Q_1^r(t)}{\mu_1^r} + \frac{Q_2^r(t)}{\mu_2^r}, \\ W_2^r(t) &\doteq \frac{Q_2^r(t)}{\mu_3^r} + \frac{Q_3^r(t)}{\mu_3^r}. \end{aligned}$$

The service allocation processes are required to satisfy the conditions below.

For $j = 1, 2, 3, k = 1, 2, r \in \mathbb{S}$,

(2.3)         $$T_j^r(t) \in \mathcal{F}; \qquad t \geq 0,$$

(2.4)     $T_j^r(\cdot)$ is a continuous nondecreasing process with $T_j^r(0) = 0$,

(2.5)     $I_k^r(\cdot)$ is a continuous nondecreasing process with $I_k^r(0) = 0$,

(2.6)         $$Q_k^r(t) \geq 0, \qquad t \geq 0.$$

From (2.3)–(2.5) and recalling the definition of $I_1^r(\cdot)$ and $I_2^r(\cdot)$, we get that, for all $j = 1, 2, 3$,

(2.7) $\quad T_j^r$ is uniformly Lipschitz continuous with Lipschitz constant 1.

Any process $T^r$ satisfying (2.3)–(2.6) will be referred to as an admissible control policy for the $r$th network. Note that we are not assuming any further measurability condition on $T^r$ except (2.3).

Now we define fluid-scaled processes and diffusion-scaled processes corresponding to the processes described above. For each $r \in \mathbb{S}$ and an admissible control policy $T^r(\cdot)$ with associated queue-length process $Q^r(\cdot)$ and idle-time process $I^r(\cdot)$, define, for $t \geq 0$,

*Fluid-scaled processes.*

$$\begin{aligned}
\bar{T}^r(t) &\doteq r^{-2} T^r(r^2 t), & \bar{I}^r(t) &\doteq r^{-2} I^r(r^2 t), \\
\bar{A}^r(t) &\doteq r^{-2} A^r(r^2 t), & \bar{S}^r(t) &\doteq r^{-2} S^r(r^2 t), \\
\bar{Q}^r(t) &\doteq r^{-2} Q^r(r^2 t), & \bar{W}^r(t) &\doteq r^{-2} W^r(r^2 t).
\end{aligned}$$

(2.8)

*Diffusion-scaled processes.*

$$\begin{aligned}
\hat{T}^r(t) &\doteq r^{-1} T^r(r^2 t), & \hat{I}^r(t) &\doteq r^{-1} I^r(r^2 t), \\
\hat{A}^r(t) &\doteq r^{-1}(A^r(r^2 t) - \lambda^r r^2 t), & \hat{S}^r(t) &\doteq r^{-1}(S^r(r^2 t) - \mu^r r^2 t), \\
\hat{Q}^r(t) &\doteq r^{-1} Q^r(r^2 t), & \hat{W}^r(t) &\doteq r^{-1} W^r(r^2 t).
\end{aligned}$$

(2.9)

By the definitions above, we have the following identities. For all $t \geq 0$,

(2.10) $$\hat{W}^r(t) = M^r \hat{Q}^r(t),$$

where

$$M^r \doteq \begin{pmatrix} \dfrac{1}{\mu_1^r} & \dfrac{1}{\mu_2^r} & 0 \\[2mm] 0 & \dfrac{1}{\mu_3^r} & \dfrac{1}{\mu_3^r} \end{pmatrix}$$

and

(2.11) $$\begin{aligned}
\hat{Q}_i^r(t) &= (\hat{A}_i^r(t) - \hat{S}_i^r(\bar{T}_i^r(t))) + r(\lambda_i^r t - \mu_i^r \bar{T}_i^r(t)), \qquad i = 1, 2, \\
\hat{Q}_3^r(t) &= (\hat{S}_2^r(\bar{T}_2^r(t)) - \hat{S}_3^r(\bar{T}_3^r(t))) + r(\mu_2^r \bar{T}_2^r(t) - \mu_3^r \bar{T}_3^r(t)).
\end{aligned}$$

We also define another process $\hat{X}^r(\cdot)$, which is closely related to the scaled queue length process $\hat{Q}^r(\cdot)$. A formal limit of $\hat{X}^r(\cdot)$ is used in the BCP described in Section 3. For $t \geq 0$, let

(2.12) $$\begin{aligned}
\hat{X}_i^r(t) &\doteq \hat{A}_i^r(t) - \hat{S}_i^r(\bar{T}_i^r(t)) + r\left(\lambda_i^r t - \mu_i^r \frac{\lambda_i}{\mu_i} t\right), \qquad i = 1, 2, \\
\hat{X}_3^r(t) &\doteq \hat{S}_2^r(\bar{T}_2^r(t)) - \hat{S}_3^r(\bar{T}_3^r(t)) + r\left(\mu_2^r \frac{\lambda_2}{\mu_2} t - \mu_3^r t\right).
\end{aligned}$$

From (2.11) and (2.12), we have the following relationships:

$$
\begin{aligned}
(2.13) \qquad \hat{Q}_i^r(t) &= \hat{X}_i^r(t) + r\mu_i^r\left(\frac{\lambda_i}{\mu_i}t - \bar{T}_i^r(t)\right), \qquad i = 1, 2, \\
\hat{Q}_3^r(t) &= \hat{X}_3^r(t) + r\mu_3^r(t - \bar{T}_3^r(t)) - r\mu_2^r\left(\frac{\lambda_2}{\mu_2}t - \bar{T}_2^r(t)\right).
\end{aligned}
$$

We will assume that the sequence of networks is in heavy traffic. More precisely, we will make the following assumption.

ASSUMPTION 2.2 (Heavy traffic assumption). We assume that the following relationships hold for the limiting parameters:

$$
(2.14) \qquad \frac{\lambda_1}{\mu_1} + \frac{\lambda_2}{\mu_2} = 1, \qquad \frac{\lambda_2}{\mu_3} = 1,
$$

and there exist $b_i \in \mathbb{R}$, $i = 1, 2, 3$, such that $\lim_{r\to\infty} b_i^r = b_i$, where

$$
(2.15) \qquad b_i^r \doteq r\left(\frac{\lambda_i^r}{\mu_i^r} - \frac{\lambda_i}{\mu_i}\right), \qquad i = 1, 2, \qquad b_3^r \doteq r\left(\frac{\lambda_2^r}{\mu_3^r} - 1\right).
$$

Under Assumptions 2.1 and 2.2, the diffusion-scaled workload process has the following representation:

$$
(2.16) \qquad \hat{W}^r(t) = M^r \hat{X}^r(t) + \hat{I}^r(t),
$$

where, for $t \geq 0$, $\hat{I}^r(t) \doteq I^r(r^2 t)/r$ and $I_1^r(t) \doteq t - T_1^r(t) - T_2^r(t), I_2^r(t) \doteq t - T_3^r(t)$.

2.4. *The cost function.* For the $r$th system, we consider the expected infinite horizon discounted (linear) holding cost associated with the control $T^r$ and the corresponding normalized queue-length process $\hat{Q}^r$, given as follows:

$$
(2.17) \qquad \hat{J}^r(T^r) = \mathbb{E}\left(\int_0^\infty e^{-\gamma t} h \cdot \hat{Q}^r(t)\, dt\right),
$$

where $\gamma \in (0, \infty)$ is the "discount factor" and $h \equiv (h_1, h_2, h_3)$; $h_k \in (0, \infty)$, $k = 1, 2, 3$, is the vector of "holding costs" for the three buffers.

The goal is to find a sequence of admissible controls which asymptotically give the minimum possible cost, that is, find a sequence $\{T^r\}$ such that

$$
\lim_{r\to\infty} \hat{J}^r(T^r) = \inf \liminf_{r\to\infty} \hat{J}^r(\tilde{T}^r),
$$

where the infimum on the right-hand side is taken over all admissible sequences $\{\tilde{T}^r\}$.

We will make the following assumption on the service rate and holding cost parameters.

ASSUMPTION 2.3.   $h_1\mu_1 - h_2\mu_2 + h_3\mu_2 > 0$, $h_2\mu_2 - h_3\mu_2 \geq 0$ and $h_2\mu_2 - h_1\mu_1 \geq 0$.

This parameter regime is the *Case IIA* of [11] among the different cases mentioned in that paper. *Case I* considers the parameter regime $h_1\mu_1 - h_2\mu_2 + h_3\mu_2 \leq 0$. This case has a simple priority policy which is shown to be asymptotically optimal in [16]. *Case II* corresponds to the complementary regime, namely $h_1\mu_1 - h_2\mu_2 + h_3\mu_2 > 0$. In this case, for the first server, serving Class 1 jobs reduces immediate cost at an average rate of $h_1\mu_1$, whereas serving Class 2 jobs would reduce immediate cost at an average rate of $h_2\mu_2$, but increases cost at an average rate of $h_3\mu_2$, since a job served from Class 2 becomes a Class 3 job. Since $h_1\mu_1 > h_2\mu_2 - h_3\mu_2$, total immediate cost is reduced at a more rapid average rate by serving Class 1 jobs. But a simple priority policy for *Server* 1 that requires it to work on Class 1 jobs, whenever *Buffer* 1 is nonempty, will cause starvation of *Server* 2 and is likely to cause the contents of *Buffer* 2 to grow without bound. In the *Case IIA*, we also assume $h_2\mu_2 \geq h_1\mu_1$ and $h_2\mu_2 \geq h_3\mu_2$ (or, simply, $h_2 \geq h_3$). Here the second condition means that it is cheaper to hold jobs in *Buffer* 3 than in *Buffer* 2. Also, the first condition above says that working on *Buffer* 2 reduces the immediate cost at *Server* 1 more quickly than working on *Buffer* 1. In this work we show that, under Assumption 2.3, a suitable threshold policy is asymptotically optimal. This policy (see Definition 3.6 for the precise description of the policy) keeps a sufficient number of jobs in *Buffer* 3 (so that *Server* 2 does not idle unnecessarily) and makes *Server* 1 work on both the associated buffers so that none of the buffers blow up. An example of parameters satisfying Assumption 2.3 is $h_1 = h_2 = h_3 = 1, \mu_1 = \mu_2 = 2, \mu_3 = 1$. In [7] the authors worked with this set of parameter values.

**3. Brownian control problem.**   We now introduce the BCP (see [3]) associated with the crisscross network introduced above. This control problem is obtained by taking a formal limit of the control problems for the above sequence of networks. More precisely, defining

$$(3.1) \qquad \bar{T}^*(t) \doteq \left(\frac{\lambda_1}{\mu_1}t, \frac{\lambda_2}{\mu_2}t, t\right), \qquad t \geq 0,$$

we might expect, for "reasonable" control policies, that, as $r \to \infty$,

$$(3.2) \qquad \bar{T}^r \Rightarrow \bar{T}^*.$$

From the functional central limit theorem, one has that

$$(3.3) \qquad (\hat{A}^r(\cdot), \hat{S}^r(\cdot)) \Rightarrow (\tilde{A}(\cdot), \tilde{S}(\cdot)),$$

where $\tilde{A}$ is a two-dimensional Brownian motion that starts from the origin and has diagonal covariance matrix, $\mathrm{diag}(\lambda_1, \lambda_2)$ and $\tilde{S}$ is a three-dimensional Brownian motion, independent of $\tilde{A}$, that starts from the origin and has diagonal covariance matrix, $\mathrm{diag}(\mu_1, \mu_2, \mu_3)$. Using (3.3), (3.2), a random time change theorem (Lemma 3.14.1 of [2]) and the heavy traffic condition (Assumption 2.2), one has that

$$\hat{X}^r(\cdot) \Rightarrow \tilde{X}(\cdot), \tag{3.4}$$

where, for $t \geq 0$,

$$\tilde{X}_i(t) = \tilde{A}_i(t) - \tilde{S}_i(T_i^*(t)) + (\mu_i b_i)t, \qquad i = 1, 2,$$
$$\tilde{X}_3(t) = \tilde{S}_2(T_2^*(t)) - \tilde{S}_3(T_3^*(t)) + (\mu_3 b_3 - \mu_2 b_2)t.$$

Note that $\tilde{X}$ is a three-dimensional Brownian motion that starts from origin, with a drift $(\mu_1 b_1, \mu_2 b_2, \mu_3 b_3 - \mu_2 b_2)$ and covariance matrix

$$\begin{pmatrix} 2\lambda_1 & 0 & 0 \\ 0 & 2\lambda_2 & -\lambda_2 \\ 0 & -\lambda_2 & 2\lambda_2 \end{pmatrix}.$$

As stated in the beginning of Section 2.2, we can assume (by enlarging the probability space, if needed), without loss of generality, that $\tilde{A}, \tilde{S}, \tilde{X}$ are defined on $(\Omega, \mathcal{F}, \mathbb{P})$. Thus, taking a formal limit as $r \to \infty$ in (2.11), (2.12), (3.2) and (2.17), one arrives at the following *Brownian control problem*.

DEFINITION 3.1 [Brownian control problem (BCP)]. Let $\tilde{X}(\cdot)$ be as defined below (3.4). The BCP is to find an $\mathbb{R}^3$-valued measurable stochastic process $\tilde{Y}(\cdot) \doteq (\tilde{Y}_1(\cdot), \tilde{Y}_2(\cdot), \tilde{Y}_3(\cdot))$, referred to as the control process, which minimizes

$$\mathbb{E}\left( \int_0^\infty e^{-\gamma t} h \cdot \tilde{Q}(t) \, dt \right), \tag{3.5}$$

subject to the following conditions. For all $t \geq 0$,

$$\begin{aligned} 0 &\leq \tilde{Q}_1(t) \doteq \tilde{X}_1(t) + \mu_1 \tilde{Y}_1(t), \\ 0 &\leq \tilde{Q}_2(t) \doteq \tilde{X}_2(t) + \mu_2 \tilde{Y}_2(t), \\ 0 &\leq \tilde{Q}_3(t) \doteq \tilde{X}_3(t) + \mu_3 \tilde{Y}_3(t) - \mu_2 \tilde{Y}_2(t) \end{aligned} \tag{3.6}$$

and

$$\begin{aligned} \tilde{I}_1(\cdot) &\doteq \tilde{Y}_1(\cdot) + \tilde{Y}_2(\cdot) \text{ is nondecreasing and } \tilde{I}_1(0) = 0, \\ \tilde{I}_2(\cdot) &\doteq \tilde{Y}_3(\cdot) \text{ is nondecreasing and } \tilde{I}_2(0) = 0. \end{aligned} \tag{3.7}$$

We will refer to any measurable process $\tilde{Y}(\cdot)$ satisfying (3.6) and (3.7) as an admissible control for the BCP.

REMARK 3.2. Our formulation of the BCP is somewhat different from that in Harrison (cf. [5]) in that we do not work with a weak formulation and we do not require the adaptedness of the control process. However, the diffusion control problem is not the real topic of interest here. It is used only to prove asymptotic optimality of our policy, namely, the result:

$$\lim_{r \to \infty} \hat{J}^r(T^r) = \inf \liminf_{r \to \infty} \hat{J}^r(\tilde{T}^r),$$

where $T^r$ is our proposed policy as in Definition 3.6 and the infimum on the right-hand side is taken over all $\tilde{T}^r$ satisfying (2.3)–(2.6). In this regard, the formulation considered in the current work suffices. It will be seen in Section 3.1 that the cost in (3.5) is minimized by $\tilde{Y}^*$ which is adapted to the filtration generated by $\tilde{X}$. It follows that the infimum of the cost in (3.5), $\tilde{J}^*$ [see (3.27)], is the same as that taken over all probability spaces supporting a three-dimensional Brownian motion with the same drift and covariance matrix as $\tilde{X}$.

3.1. *Reduction to the equivalent workload formulation.* Let $\tilde{Y}(\cdot)$ be an admissible control for the BCP and define $\tilde{Q}$ via (3.6). Define the workload process

$$(3.8) \qquad \tilde{W}(t) \doteq M\tilde{Q}(t), \qquad t \geq 0,$$

where

$$M \doteq \begin{pmatrix} \dfrac{1}{\mu_1} & \dfrac{1}{\mu_2} & 0 \\[2mm] 0 & \dfrac{1}{\mu_3} & \dfrac{1}{\mu_3} \end{pmatrix}.$$

Thus, for $t \geq 0$,

$$(3.9) \qquad \begin{aligned} \tilde{W}_1(t) &= \frac{\tilde{Q}_1(t)}{\mu_1} + \frac{\tilde{Q}_2(t)}{\mu_2}, \\ \tilde{W}_2(t) &= \frac{\tilde{Q}_2(t)}{\mu_3} + \frac{\tilde{Q}_3(t)}{\mu_3}. \end{aligned}$$

It is easy to check that

$$(3.10) \qquad \tilde{W} = M\tilde{X} + \tilde{V} \qquad \text{with } \tilde{V}(t) = (\tilde{I}_1(t), \tilde{I}_2(t))', t \geq 0.$$

We will now obtain a solution of the BCP using the above workload process. We begin by considering the following simple linear programming problem. Fix $w_1, w_2 \in [0, \infty)$. The linear program (LP) problem is as follows:

$$(3.11) \qquad \begin{aligned} &\text{minimize}_{z_1, z_2, z_3} \quad h_1 z_1 + h_2 z_2 + h_3 z_3 \\ &\text{subject to} \quad \frac{z_1}{\mu_1} + \frac{z_2}{\mu_2} = w_1, \\ &\hspace{4.5em} \frac{z_2}{\mu_3} + \frac{z_3}{\mu_3} = w_2, \\ &\hspace{4.5em} z_1, z_2, z_3 \geq 0. \end{aligned}$$

A straightforward calculation using the fact that $h_1\mu_1 - h_2\mu_2 + h_3\mu_2 > 0$ (see Assumption 2.3) shows (cf. [11]) that the value of the LP is

(3.12)    $$\hat{h}(w_1, w_2) = \begin{cases} (h_2\mu_2 - h_3\mu_2)w_1 + (h_3\mu_3)w_2, \\ \qquad \text{when } \mu_3 w_2 \geq \mu_2 w_1, \\ (h_1\mu_1)w_1 + \dfrac{\mu_3}{\mu_2}(h_2\mu_2 - h_1\mu_1)w_2, \\ \qquad \text{when } \mu_3 w_2 \leq \mu_2 w_1. \end{cases}$$

In particular, if $z_1, z_2, z_3$ are nonnegative numbers such that $\frac{z_1}{\mu_1} + \frac{z_2}{\mu_2} = w_1$ and $\frac{z_2}{\mu_3} + \frac{z_3}{\mu_3} = w_2$, then

(3.13)    $$h_1 z_1 + h_2 z_2 + h_3 z_3 \geq \hat{h}(w_1, w_2).$$

Another simple calculation yields the following solution of the LP:

(3.14)
$$\begin{aligned} &z_1^* = 0, &&z_2^* = \mu_2 w_1, &&z_3^* = \mu_3 w_2 - \mu_2 w_1, \\ &&&&&\text{if } \mu_3 w_2 \geq \mu_2 w_1, \\ &z_1^* = \frac{\mu_1}{\mu_2}(\mu_2 w_1 - \mu_3 w_2), &&z_2^* = \mu_3 w_2, &&z_3^* = 0, \\ &&&&&\text{if } \mu_3 w_2 \leq \mu_2 w_1. \end{aligned}$$

REMARK 3.3.  Note that from Assumption 2.3, $h_2\mu_2 - h_3\mu_2 \geq 0$ and $h_2\mu_2 - h_1\mu_1 \geq 0$. Thus, we have that $\hat{h}(w_1, w_2)$ is a nondecreasing function of both $w_1$ and $w_2$. This monotonicity property is critical in obtaining a pathwise optimal solution to the BCP.

We now present another control problem which, because of the monotonicity property of $\hat{h}$, can be solved explicitly. The results of [6] show that, using a solution of this reduced control problem (referred to as EWF in Definition 3.4 below) and the solution $\hat{h}$ of the linear program in (3.11), one can obtain a solution of the BCP.

DEFINITION 3.4 [Equivalent workload formulation (EWF)].  Let $\tilde{X}(\cdot)$ be as defined below (3.4). The equivalent workload problem is to find an $\mathbb{R}^2$-valued measurable stochastic process $\tilde{I}(\cdot) = (\tilde{I}_1(\cdot), \tilde{I}_2(\cdot))$, referred to as the control process, which minimizes

(3.15)    $$\mathbb{E}\left( \int_0^\infty e^{-\gamma t} \hat{h}(\tilde{W}(t))\, dt \right),$$

subject to the following conditions. For all $t \geq 0$

(3.16)
$$\begin{aligned} &0 \leq \tilde{W}_1(t) \doteq \frac{\tilde{X}_1(t)}{\mu_1} + \frac{\tilde{X}_2(t)}{\mu_2} + \tilde{I}_1(t), \\ &0 \leq \tilde{W}_2(t) \doteq \frac{\tilde{X}_2(t)}{\mu_3} + \frac{\tilde{X}_3(t)}{\mu_3} + \tilde{I}_2(t) \quad \text{and} \\ &\tilde{I}_1(\cdot), \tilde{I}_2(\cdot) \text{ are nondecreasing and } \tilde{I}(0) = 0. \end{aligned}$$

From (3.16), using the minimality property of the one-dimensional Skorohod problem (see Proposition B.1 of [1]), we have that if

$$(3.17) \qquad (\tilde{V}_1^*(t), \tilde{V}_2^*(t)) \doteq \left( -\inf_{0 \leq s \leq t}(m_1' \tilde{X}(s)), -\inf_{0 \leq s \leq t}(m_2' \tilde{X}(s)) \right)',$$

then

$$(3.18) \qquad \tilde{W}_i(t) \geq \tilde{W}_i^*(t) \doteq m_i' \tilde{X}(t) + \tilde{V}_i^*(t), \qquad i = 1, 2,$$

where, for $i = 1, 2$, $m_i'$ is the $i$th row of the matrix $M$, that is, $M = [m_1 : m_2]'$. Also, from (3.18) and Remark 3.3, it follows that, for all $t \geq 0$,

$$(3.19) \qquad \hat{h}(\tilde{W}_1(t), \tilde{W}_2(t)) \geq \hat{h}(\tilde{W}_1^*(t), \tilde{W}_2^*(t)).$$

This shows that $\tilde{V}^*$ is a solution to the EWF. Using $\tilde{W}^*$ and $\tilde{V}^*$, we now construct the solution of the BCP. The solution is motivated by the solution of the LP problem in (3.11), given via (3.14). Define processes $\tilde{Y}_i^*(\cdot)$, $i = 1, 2, 3$, as follows. For $t \geq 0$, let

$$(3.20) \quad \tilde{Y}_1^*(t) \doteq \begin{cases} -\dfrac{\tilde{X}_1(t)}{\mu_1}, & \text{if } \mu_3 \tilde{W}_2^*(t) \geq \mu_2 \tilde{W}_1^*(t), \\ -\dfrac{\tilde{X}_3(t)}{\mu_2} + \tilde{V}_1^*(t) - \dfrac{\mu_3}{\mu_2} \tilde{V}_2^*(t), & \text{if } \mu_3 \tilde{W}_2^*(t) < \mu_2 \tilde{W}_1^*(t), \end{cases}$$

$$(3.21) \quad \tilde{Y}_2^*(t) \doteq \begin{cases} \dfrac{\tilde{X}_1(t)}{\mu_1} + \tilde{V}_1^*(t), & \text{if } \mu_3 \tilde{W}_2^*(t) \geq \mu_2 \tilde{W}_1^*(t), \\ \dfrac{\tilde{X}_3(t)}{\mu_2} + \dfrac{\mu_3}{\mu_2} \tilde{V}_2^*(t), & \text{if } \mu_3 \tilde{W}_2^*(t) < \mu_2 \tilde{W}_1^*(t), \end{cases}$$

and

$$(3.22) \qquad \tilde{Y}_3^*(t) \doteq \tilde{V}_2^*(t).$$

It is easy to verify that $\tilde{Y}^*$ is an admissible control for the BCP. Also, it follows from (3.20)–(3.22) and (3.16) that $\tilde{I}^* = \tilde{V}^*$. Now define $\tilde{Q}^*$ via (3.6), with $\tilde{Y}$ there replaced by $\tilde{Y}^*$.

Hence, we have that if $\mu_3 \tilde{W}_2^*(t) \geq \mu_2 \tilde{W}_1^*(t)$, then

$$(3.23) \quad \begin{aligned} \tilde{Q}_1^*(t) &= 0, \\ \tilde{Q}_2^*(t) &= \mu_2 \tilde{W}_1^*(t) \\ &= \left( \frac{\mu_2}{\mu_1} \right) \tilde{X}_1(t) + \tilde{X}_2(t) + \mu_2 \tilde{V}_1^*(t), \\ \tilde{Q}_3^*(t) &= \mu_3 \tilde{W}_2^*(t) - \mu_2 \tilde{W}_1^*(t) \\ &= -\left( \frac{\mu_2}{\mu_1} \right) \tilde{X}_1(t) + \tilde{X}_3(t) + \mu_3 \tilde{V}_2^*(t) - \mu_2 \tilde{V}_1^*(t), \end{aligned}$$

and if $\mu_3 \tilde{W}_2^*(t) < \mu_2 \tilde{W}_1^*(t)$, then

$$
\begin{aligned}
(3.24) \quad \tilde{Q}_1^*(t) &= \left(\frac{\mu_1}{\mu_2}\right)(\mu_2 \tilde{W}_1^*(t) - \mu_3 \tilde{W}_2^*(t)) \\
&= \tilde{X}_1(t) - \left(\frac{\mu_1}{\mu_2}\right)\tilde{X}_3(t) + \mu_1 \tilde{V}_1^*(t) - \frac{\mu_1\mu_3}{\mu_2}\tilde{V}_2^*(t), \\
\tilde{Q}_2^*(t) &= \mu_3 \tilde{W}_2^*(t) \\
&= (\tilde{X}_2(t) + \tilde{X}_3(t)) + \mu_3 \tilde{V}_2^*(t), \\
\tilde{Q}_3^*(t) &= 0.
\end{aligned}
$$

Now we show that $\tilde{Y}^*$ is a solution to the BCP described in the beginning of this section. Note that if $\tilde{Y}$ is any admissible control for the BCP and $\tilde{Q}$ is defined via (3.6), then, from (3.13), for all $t \geq 0$,

$$
(3.25) \qquad h \cdot \tilde{Q}(t) \geq \hat{h}(\tilde{W}_1(t), \tilde{W}_2(t)),
$$

where $\tilde{W}_1$ and $\tilde{W}_2$ are defined via (3.9). In view of (3.25) and (3.19), in order to show that $\tilde{Y}^*$ is the solution of the BCP, it suffices to show that

$$
(3.26) \qquad \hat{h}(\tilde{W}_1^*(t), \tilde{W}_2^*(t)) = h \cdot \tilde{Q}^*(t).
$$

But (3.26) is an immediate consequence of the definition of $\tilde{Q}^*$ [see (3.23) and (3.24) and the fact that $z^*$ defined via (3.14) is the solution of the linear program in (3.11)].

This proves that $\tilde{Y}^*$ is a solution for the BCP with the corresponding queue-length $\tilde{Q}^*$. Let the infimum of the objective function (3.5), over all admissible controls, in the BCP be denoted by $J^*$, that is,

$$
(3.27) \qquad J^* \doteq \inf \mathbb{E}\left(\int_0^\infty e^{-\gamma t} h \cdot \tilde{Q}(t)\,dt\right).
$$

Thus, we have that

$$
(3.28) \qquad J^* = \mathbb{E}\left(\int_0^\infty e^{-\gamma t} h \cdot \tilde{Q}^*(t)\,dt\right).
$$

REMARK 3.5. The BCP (Definition 3.1) presented in this work is somewhat different from the BCP studied in Section 7 of [11]. The BCP considered in [11] is formulated in terms of a four-dimensional control process which is required to have paths of bounded variation. Due to this restriction, the authors were unable to prove the existence of an optimal control policy for the BCP. For precise description of the control problem, we refer the readers to Section 7 of [11]. In the formulation considered in the current paper, the BCP has a three-dimensional control process which is restricted to have paths in $\mathcal{D}([0,\infty); \mathbb{R}^3)$. Thus is our formulation, an admissible control need not have paths of bounded variation. As seen above, the BCP in Definition 3.1 has an optimal solution given via (3.20)–(3.22).

3.2. *The policy.* Motivated by the solution of the BCP, we now propose our control policy for the $r$th network, $r \in \mathbb{S}$. Fix $c, \ell_0 \in (1, \infty)$. Define $\mathbf{L^r} \doteq \lfloor \ell_0 \log r \rfloor$ and $\mathbf{C^r} \doteq \lfloor c_0 \log r \rfloor$, where $c_0 = c\ell_0$. Since we are interested in asymptotic optimality, we can (and will) assume, without loss of generality, that $r \geq \bar{r}$, where $\bar{r}$ is such that for all $r \geq \bar{r}$, $\mathbf{C^r} - \mathbf{L^r} - 1 \geq 1$ and $\frac{\mu_1^r}{\mu_2^r}(\mathbf{C^r} - \mathbf{L^r} + 2) \geq 1$.

DEFINITION 3.6 (Control policy). The policy is as follows. No idling by *Server* 2 unless *Buffer* 3 is empty. The sequencing control for *Server* 1 is as follows.

If $Q_3^r(s) - \frac{\mu_2^r}{\mu_1^r}Q_1^r(s) < \mathbf{L^r}$,

serve *Buffer* 2 if $Q_3^r(s) < \mathbf{C^r} - 1$ and $Q_2^r(s) \neq 0$,

serve *Buffer* 1 (when it is nonempty) if either $Q_3^r(s) \geq \mathbf{C^r} - 1$ or $Q_2^r(s) = 0$.

If $Q_3^r(s) - \frac{\mu_2^r}{\mu_1^r}Q_1^r(s) \geq \mathbf{L^r}$,

serve *Buffer* 1 (when it is nonempty) if either $Q_1^r(s) \geq \frac{\mu_1^r}{\mu_2^r}(\mathbf{C^r} - \mathbf{L^r} + 2)$ or $Q_2^r(s) = 0$, serve *Buffer* 2 if $Q_1^r(s) < \frac{\mu_1^r}{\mu_2^r}(\mathbf{C^r} - \mathbf{L^r} + 2)$ and $Q_2^r(s) \neq 0$.

*Server* 1 idles if both *Buffer* 1 and *Buffer* 2 are empty.

We will refer to the constants $c$ and $\ell_0$ as the threshold parameters of the control policy. It will be shown that, for a choice of $c$ and $\ell_0$ large enough, the above policy is asymptotically optimal. One precise choice of $c$ and $\ell_0$ is given in Remark 4.3(a).

One of the referees has conjectured that the above policy with $\ell_0 = 0$ and $c_0$ replaced by a sufficiently large constant is asymptotically optimal as well. However, as is explained in the following paragraph, the arguments in the current paper crucially rely on the largeness of $\ell_0$.

Now we provide some motivation for the policy proposed above. Note that

$$\mu_3^r \hat{W}_2^r(t) - \mu_2^r \hat{W}_1^r(t) = \hat{Q}_3^r(t) - \frac{\mu_2^r}{\mu_1^r}\hat{Q}_1^r(t).$$

Thus, the solution of the BCP suggests that when $\hat{Q}_3^r(t) - \frac{\mu_2^r}{\mu_1^r}\hat{Q}_1^r(t) < 0$, then the optimal policy should try to make *queue* 3 empty, whereas when the opposite is true, *queue* 1 should be emptied. This is achieved in the first regime via the threshold $\mathbf{C^r} - 1$ and in the other regime via the threshold $\frac{\mu_1^r}{\mu_2^r}(\mathbf{C^r} - \mathbf{L^r} + 2)$. Note that the two thresholds $\mathbf{C^r} - 1$ and $\frac{\mu_1^r}{\mu_2^r}(\mathbf{C^r} - \mathbf{L^r} + 2)$ approach $\infty$ as $r \to \infty$, however, in diffusion scaling these are negligible. Furthermore, (3.19), suggests that asymptotically there should be no idling by *Server* 1 unless there is no work in *Buffer* 1 and *Buffer* 2, and that there be no idling by *Server* 2 unless there is no work at *Buffer* 2 and *Buffer* 3. The first nonidleness condition is quite easy to enforce, by saying that the first server works whenever there is work for it to do. However, the second nonidleness condition is difficult to enforce, since one can get into the situation where *Buffer* 3 is empty and so *Server* 2 has no immediate work to do but *Buffer* 2 is nonempty. Thus, one needs to ensure that there is always enough work in *Buffer* 3 when *Buffer* 2 is nonempty. This is the reason for the threshold $\mathbf{L^r} = \lfloor \ell_0 \log r \rfloor$ in the policy. For our proof of asymptotic optimality, we will need that $\ell_0$ is sufficiently large (see Theorem 4.9).

REMARK 3.7. This policy is preemptive-resume type. For example, if at any time instant $t$ *Server* 1 is working on jobs of Class 1 and the policy requires it to work on Class 2 jobs, it immediately suspends all Class 1 jobs and starts working on Class 2 jobs (suspended jobs if there are any, or new jobs). When at a later time it turns to Class 1 jobs again, it resumes working on the suspended Class 1 job (and spends only the excess time that it needs to complete the remaining part of the job, so that the total time spent on this job is the same as the time needed to complete this job if there was no interruption).

REMARK 3.8. The above policy can be written in the following form. Let

$$
\begin{aligned}
\mathcal{A}_r &= \left\{ u \geq 0 : Q_3^r(u) - \frac{\mu_2^r}{\mu_1^r} Q_1^r(u) < \mathbf{L^r} \right\}, \\
\mathcal{B}_r &= \{ u \geq 0 : Q_3^r(u) \geq \mathbf{C^r} - 1 \text{ or } Q_2^r(u) = 0 \}, \\
\mathcal{C}_r &= \left\{ u \geq 0 : Q_1^r(u) \geq \frac{\mu_1^r}{\mu_2^r} (\mathbf{C^r} - \mathbf{L^r} + 2) \text{ or } Q_2^r(u) = 0 \right\}, \\
\mathcal{D}_r &= \{ u \geq 0 : Q_1^r(u) + Q_2^r(u) \neq 0 \}.
\end{aligned}
\tag{3.29}
$$

Then the (sequence of) proposed policies $\{T^r\}$ described above in Definition 3.6 can be described as follows. For $j = 1, 2, 3$, $T_j^r$ is the absolutely continuous function whose derivative (defined a.e.), denoted by $\dot{T}_j^r$, is given as follows:

$$
\begin{aligned}
\dot{T}_1^r &= (\mathbf{I}_{\mathcal{A}_r} \mathbf{I}_{\mathcal{B}_r} + \mathbf{I}_{\mathcal{A}_r^c} \mathbf{I}_{\mathcal{C}_r}) \mathbf{I}_{\mathcal{D}_r}, \\
\dot{T}_2^r &= (\mathbf{I}_{\mathcal{A}_r} \mathbf{I}_{\mathcal{B}_r^c} + \mathbf{I}_{\mathcal{A}_r^c} \mathbf{I}_{\mathcal{C}_r^c}) \mathbf{I}_{\mathcal{D}_r}, \\
\dot{T}_3^r &= \mathbf{I}_{\{u \geq 0 : Q_5^r(u) > 0\}}.
\end{aligned}
\tag{3.30}
$$

Note that, for $j = 1, 2, 3$,

$$
T_j^r(t) = \int_0^t \dot{T}_j^r(s) \, ds.
$$

Here $\mathbf{I}_A$ denotes the indicator function of a set $A$ and $A^{\mathbf{c}}$ denotes the complement of a set $A$. Note that $\dot{T}_3^r(t)$ and $\dot{T}_1^r(t) + \dot{T}_2^r(t)$ are both $\{0, 1\}$ valued and $\dot{T}_1^r(t) + \dot{T}_2^r(t) = 0$ if and only if both $Q_1^r(t)$ and $Q_2^r(t)$ are zero, and $\dot{T}_3^r(t) = 0$ if and only if $Q_3^r(t)$ is zero. In other words, the policy operates in a "nonidling" fashion.

**4. Proof of asymptotic optimality of the proposed policy.** In this section we will prove the asymptotic optimality of the scheduling control policy introduced in Definition 3.6. More precisely, we prove the following two results.

THEOREM 4.1. *Let $\{T^r\}$ be any sequence of scheduling controls. Then, for $J^*$ as in (3.28), we have*

$$\liminf_{r \to \infty} \hat{J}^r(T^r) \geq J^*. \tag{4.1}$$

THEOREM 4.2. *There exist $c, \bar{\ell} \in (1, \infty)$ such that if $\{T^r\}$ is the sequence of scheduling controls described in Definition 3.6 with threshold parameters $c$ and $\ell_0$ with $\ell_0 \geq \bar{\ell}$, then*

$$\begin{array}{ll} \text{(a)} & (\hat{W}^r, \hat{I}^r) \Rightarrow (\tilde{W}^*, \tilde{I}^*) \text{ as } r \to \infty, \\ \text{(b)} & \lim_{r \to \infty} \hat{J}^r(T^r) = J^*, \end{array} \tag{4.2}$$

*where $J^*$ is as in (3.28).*

Theorem 4.1 says that the asymptotic cost associated with any scheduling policy cannot be lower than $J^*$ defined in (3.28), while Theorem 4.2 says that the control described in Definition 3.6 asymptotically achieves $J^*$, which is the optimal cost for the BCP.

REMARK 4.3. (a) The choice of $c, \bar{\ell}$ depends on various large deviation estimates that are obtained in Section 5. A concrete choice of $c, \bar{\ell} \in (1, \infty)$ is as follows:

$$\begin{array}{l} \theta_3 = \dfrac{\mu_1}{\mu_2 \lambda_1} \min\{\eta_1, \eta_2\}, \\ \rho_2 = \min\{\eta_3, \eta_4\}, \end{array} \tag{4.3}$$

where $\eta_1 \doteq \varsigma_2(\{\lambda_1^r\}, 1/2)$, $\eta_2 \doteq \varsigma_2(\{\mu_1^r\}, 1/2)$, $\eta_3 \doteq \varsigma_2(\{\mu_3^r\}, \min\{\mu_3/2, 1\})$, $\eta_4 \doteq \varsigma_2(\{\mu_2^r\}, \min\{\mu_2/2, 1\})$ and $\varsigma_2(\cdot)$ is as in Corollary 5.3. Choose $c, K, d, \theta$ as in (5.33). Define $\gamma_4 = (2d/K)\theta\rho_2$ and choose $\bar{\ell} = \max\{4/\gamma_4, 4/(\theta_3(c-1))\} + 1$.

(b) In this paper we restrict ourselves to a discounted cost, however, similar results can be proved for some other cost criterion (with linear holding cost) as well, by suitable modifications. The key obstacle is to prove the uniform integrability estimates of Section 5. In particular, if the criterion is finite time horizon total cost, then the uniform integrability estimates are easy to obtain.

(c) In this work we will also establish (see Corollary 4.10) that, under the proposed policy (Definition 3.6),

$$\hat{Q}_1^r(\cdot)\hat{Q}_3^r(\cdot) \Rightarrow 0.$$

In [11] the authors conjectured that any optimal policy should try to get the queue-lengths close to the set

$$\{(z_1, z_2, z_3) \in [0, \infty)^3 : z_1 z_3 = 0\}.$$

*Outline of the proofs.* The main steps in the proof of asymptotic optimality of the proposed policy are as follows. As a first step we show in Theorem 4.1 that the asymptotic cost for any sequence of policies is bounded below by $J^*$. The key step is proving the inequality in (4.12) and the main ingredients in its proof are the monotonicity property described in Remark 3.3 and the minimality property of the Skorohod map [see (4.22)]. We next show that the asymptotic cost for the sequence $\{T^{r,*}\}$ is $J^*$. The first step in this direction is obtaining the following convergence results for the queue-length and idle-time processes (see Corollary 4.10):

$$\hat{Q}_1^r(\cdot)\mathbf{I}_{\{\hat{Q}_3^r(\cdot)-(\mu_2^r/\mu_1^r)\hat{Q}_1^r(\cdot)\geq \mathbf{L}^r/r\}} \Rightarrow 0,$$

$$\hat{Q}_3^r(\cdot)\mathbf{I}_{\{\hat{Q}_3^r(\cdot)-(\mu_2^r/\mu_1^r)\hat{Q}_1^r(\cdot)< \mathbf{L}^r/r\}} \Rightarrow 0,$$

$$\int_{[0,\cdot)} \mathbf{I}_{\{\hat{Q}_2^r(s)\geq d\ell_0\log r/r\}}\, d\hat{I}_2^r(s) \Rightarrow 0.$$

The first two above are consequences of Theorem 4.8, while the third convergence result follows from Theorem 4.9. The latter result, along with the continuity of the Skorohod map, is then used to show that $(\hat{W}^r, \hat{I}^r) \Rightarrow (\tilde{W}^*, \tilde{I}^*)$. We are unable to conclude from the above convergence that $\hat{Q}^r \Rightarrow \tilde{Q}^*$; the main obstacle is showing that

$$\hat{W}_1^r(\cdot)\mathbf{I}_{\{\mu_3^r\hat{W}_3^r(\cdot)-\mu_2^r\hat{W}_1^r(\cdot)\geq \mathbf{L}^r/r\}} \Rightarrow \tilde{W}_1^*(\cdot)\mathbf{I}_{\{\mu_3\tilde{W}_3^*(\cdot)-\mu_2\tilde{W}_1^*(\cdot)\geq 0\}}.$$

However, using an elementary lemma (Lemma 4.7), we show that the convergence in (4.63) holds. Since we are working with an expected cost criterion with an unbounded cost function, in addition to the above weak convergence results, we also need suitable uniform integrability estimates. These estimates are obtained in Theorem 4.11. As an immediate consequence we then have (4.64) and (4.66). Combining these, we obtain (4.67). This along with the first two convergence results in Corollary 4.10 and the uniform integrability estimates yield (4.70). The convergence of $\hat{J}^{r,*}$ to $J^*$ then follows readily.

We begin with the following definition. Let $\mathcal{C}^m$ be the space of continuous functions from $[0, \infty)$ to $\mathbb{R}^m$ with the usual topology of uniform convergence on compact time intervals. We will suppress $m$ from the notation unless necessary.

DEFINITION 4.4 ($\mathcal{C}$-tightness).  A sequence of processes with paths in $\mathcal{D}^m$ ($m \geq 1$) is called $\mathcal{C}$-tight if it is tight in $\mathcal{D}^m$ and any weak limit point of the sequence has paths in $\mathcal{C}^m$ almost surely.

The following two basic lemmas are important in proving the optimality of the proposed policy. The proofs of the these results are similar to the proofs of Lemma 9.2 and Lemma 9.3 of [1]. However, for the sake of completeness, we have included the proofs in the Appendix.

LEMMA 4.5.  *Let $\{T^r\}$ be any sequence of scheduling policies. Then*

$$(4.4) \qquad \{\bar{Q}^r(\cdot), \bar{A}^r(\cdot), \bar{S}^r(\cdot), \bar{T}^r(\cdot), \bar{I}^r(\cdot)\}_{r \in \mathbb{S}} \ \text{is } \mathcal{C}\text{-tight}.$$

LEMMA 4.6.  *Let $\{T^r\}$ be any sequence of scheduling policies with the following property:*

$$(4.5) \qquad \underline{J}(\{T^r\}) \doteq \liminf_{r \to \infty} \hat{J}^r(T^r) < \infty,$$

*where $\hat{J}^r(T^r)$ is as in* (2.17). *Consider a subsequence $\{T^{r'}\}$ of $\{T^r\}$ such that*

$$(4.6) \qquad \lim_{r' \to \infty} \hat{J}^{r'}(T^{r'}) = \underline{J}(\{T^r\}).$$

*Then we have*

$$(4.7) \qquad \begin{aligned} (\bar{Q}^{r'}(\cdot), \bar{A}^{r'}(\cdot), \bar{S}^{r'}(\cdot), \bar{T}^{r'}(\cdot), \bar{I}^{r'}(\cdot)) \\ \Rightarrow (\mathbf{0}, \lambda(\cdot), \mu(\cdot), \bar{T}^*(\cdot), \mathbf{0}) \qquad \text{as } r' \to \infty, \end{aligned}$$

*where $\bar{T}^*$ is as defined in* (3.1), *$\mathbf{0}$ is the constant process that is zero for all $t \geq 0$, $\lambda(t) \doteq \lambda t, \mu(t) \doteq \mu t, t \geq 0$.*

Now we are ready to prove Theorem 4.1.

PROOF OF THEOREM 4.1.  If $\liminf_{r \to \infty} \hat{J}^r(T^r) = \infty$, then (4.1) holds trivially and so we only consider the case when $\liminf_{r \to \infty} \hat{J}^r(T^r) < \infty$.
Consider a subsequence $\{T^{r'}\}$ of $\{T^r\}$ such that

$$(4.8) \qquad \lim_{r' \to \infty} \hat{J}^{r'}(T^{r'}) = \liminf_{r \to \infty} \hat{J}^r(T^r) < \infty.$$

By Lemma 4.6 and (3.3), we have that, as $r' \to \infty$,

$$(4.9) \qquad (\hat{A}^{r'}(\cdot), \hat{S}^{r'}(\cdot), \bar{T}^{r'}(\cdot)) \Rightarrow (\tilde{A}(\cdot), \tilde{S}(\cdot), \bar{T}^*(\cdot)).$$

Using this observation along with Lemma 3.14.1 of [2] and Assumption 2.2 in (2.12), we have that

$$(\hat{A}^{r'}(\cdot), \hat{S}^{r'}(\cdot), \bar{T}^{r'}(\cdot), \hat{X}^{r'}(\cdot)) \Rightarrow (\tilde{A}(\cdot), \tilde{S}(\cdot), \bar{T}^*(\cdot), \tilde{X}(\cdot)),$$

where $\tilde{X}(\cdot)$ is as defined below (3.4). Using the Skorohod representation theorem, we can assume, without loss of generality, that, as $r' \to \infty$,

$$(4.10) \quad (\hat{A}^{r'}(\cdot), \hat{S}^{r'}(\cdot), \bar{T}^{r'}(\cdot), \hat{X}^{r'}(\cdot)) \to (\tilde{A}(\cdot), \tilde{S}(\cdot), \bar{T}^*(\cdot), \tilde{X}(\cdot)) \qquad \text{a.s.,}$$

uniformly on compacts (u.o.c.).

From the definition of the cost function $\hat{J}^r$ given in (2.17) and Fatou's lemma, we get

$$(4.11) \qquad \lim_{r' \to \infty} \hat{J}^{r'}(T^{r'}) \geq \mathbb{E}\left( \int_0^\infty e^{-\gamma t} \liminf_{r' \to \infty} (h \cdot \hat{Q}^{r'}(t)) \, dt \right).$$

Thus, in order to prove (4.1), it suffices to show that, for a.e. $\omega \in \Omega$ and all $t \geq 0$,

$$(4.12) \qquad \liminf_{r' \to \infty} (h \cdot \hat{Q}^{r'}(t, \omega)) \geq h \cdot \tilde{Q}^*(t, \omega),$$

where $\tilde{Q}^*(t)$ are given via the formulae in (3.17), (3.23) and (3.24) in terms of $\tilde{X}(\cdot)$ in (4.10). Fix $\omega \in \Omega$ such that $\omega$ is in the set of probability 1 on which the u.o.c. convergence in (4.10) hold, and fix $t \geq 0$. Consider the following two cases:

$$\text{Case I: } \mu_3 \tilde{W}_2^*(t, \omega) \geq \mu_2 \tilde{W}_1^*(t, \omega), \qquad \text{Case II: } \mu_3 \tilde{W}_2^*(t, \omega) < \mu_2 \tilde{W}_1^*(t, \omega),$$

where $\tilde{W}_i^*(\cdot), i = 1, 2$, are defined in terms of $\tilde{X}$ in (4.10) via the relations (3.17) and (3.18). Note that, since we are invoking the Skorohod representation theorem in (4.10), this $\tilde{W}^*$ is not the same process as in (3.18), but it has the same law as $\tilde{W}^*$ in (3.18). Once again, we retain the same symbol in order to simplify the notation.

Define $h^{r,1} \equiv (h_1^{r,1}, h_2^{r,1}, h_3^{r,1})$ as follows:

$$(4.13) \quad h_1^{r,1} \doteq \frac{h_1 \mu_1}{\mu_1^r}, \qquad h_2^{r,1} \doteq \frac{h_3 \mu_3}{\mu_3^r} + \frac{\mu_2(h_2 - h_3)}{\mu_2^r}, \qquad h_3^{r,1} \doteq \frac{h_3 \mu_3}{\mu_3^r}.$$

Observe that by Assumption 2.1, as $r \to \infty$,

$$(4.14) \qquad h_i^{r,1} \to h_i, \qquad i = 1, 2, 3.$$

From the definition of $h^{r,1}$, Assumption 2.3 and (2.10), we get

$$
\begin{aligned}
(4.15) \quad h^{r,1} \cdot \hat{Q}^r(t, \omega) &= h_1^{r,1} \hat{Q}_1^r(t, \omega) + h_2^{r,1} \hat{Q}_2^r(t, \omega) + h_3^{r,1} \hat{Q}_3^r(t, \omega) \\
&= (h_1 \mu_1)\left[ \frac{\hat{Q}_1^r(t, \omega)}{\mu_1^r} \right] + \mu_2(h_2 - h_3)\left[ \frac{\hat{Q}_2^r(t, \omega)}{\mu_2^r} \right] \\
&\quad + (h_3 \mu_3)\left[ \frac{\hat{Q}_2^r(t, \omega)}{\mu_3^r} + \frac{\hat{Q}_3^r(t, \omega)}{\mu_3^r} \right] \\
&\geq [\mu_2(h_2 - h_3)]\hat{W}_1^r(t, \omega) + [h_3 \mu_3]\hat{W}_2^r(t, \omega) \\
&= a_1^* \hat{W}_1^r(t, \omega) + b_1^* \hat{W}_2^r(t, \omega),
\end{aligned}
$$

where $a_1^* \doteq \mu_2(h_2 - h_3)$ and $b_1^* \doteq h_3\mu_3$.

Next, define $h^{r,2} \equiv (h_1^{r,2}, h_2^{r,2}, h_3^{r,2})$ as follows:

$$(4.16) \quad h_1^{r,2} = \frac{h_1\mu_1}{\mu_1^r}, \qquad h_2^{r,2} = \frac{h_1\mu_1}{\mu_2^r} + \frac{\mu_3}{\mu_3^r}\left(h_2 - \frac{h_1\mu_1}{\mu_2}\right), \qquad h_3^{r,2} = \frac{h_3\mu_3}{\mu_3^r}.$$

Once more, by Assumption 2.1, as $r \to \infty$,

$$(4.17) \qquad\qquad h_i^{r,2} \to h_i, \qquad i = 1, 2, 3,$$

and from the definition of $h^{r,2}$, Assumption 2.3 and (2.10), we get

$$
\begin{aligned}
(4.18) \quad h^{r,2} \cdot \hat{Q}^r(t,\omega) &= h_1^{r,2}\hat{Q}_1^r(t,\omega) + h_2^{r,2}\hat{Q}_2^r(t,\omega) + h_3^{r,2}\hat{Q}_3^r(t,\omega) \\
&= (h_1\mu_1)\left[\frac{\hat{Q}_1^r(t,\omega)}{\mu_1^r} + \frac{\hat{Q}_2^r(t,\omega)}{\mu_2^r}\right] \\
&\quad + \frac{\mu_3}{\mu_2\mu_3^r}[(h_2\mu_2 - h_1\mu_1)\hat{Q}_2^r(t,\omega) + (h_3\mu_3)\hat{Q}_3^r(t,\omega)] \\
&\geq (h_1\mu_1)\hat{W}_1^r(t,\omega) + \left[\frac{\mu_3}{\mu_2}(h_2\mu_2 - h_1\mu_1)\right]\hat{W}_2^r(t,\omega) \\
&= a_2^*\hat{W}_1^r(t,\omega) + b_2^*\hat{W}_2^r(t,\omega),
\end{aligned}
$$

where $a_2^* \doteq h_1\mu_1$ and $b_2^* \doteq \mu_3(h_2\mu_2 - h_1\mu_1)/\mu_2$.

Thus, defining

$$k(\omega,t) \doteq \begin{cases} 1, & \text{if } \mu_3\bar{W}_2^*(t,\omega) \geq \mu_2\bar{W}_1^*(t,\omega), \\ 2, & \text{if } \mu_3\bar{W}_2^*(t,\omega) < \mu_2\bar{W}_1^*(t,\omega), \end{cases}$$

we have that

$$(4.19) \qquad h^{r,k(\omega,t)} \cdot \hat{Q}^r(t,\omega) \geq a_{k(\omega,t)}^*\hat{W}_1^r(t,\omega) + b_{k(\omega,t)}^*\hat{W}_2^r(t,\omega).$$

Note that, from Assumption 2.3, we have that $a_i^*, b_i^*$ are nonnegative. Since $\hat{W}^r(t)$ is nonnegative for all $t \geq 0$ and $\hat{I}_1^r, \hat{I}_2^r$ are nondecreasing and start from zero, we have from (2.16) and the minimality of the solution of the Skorohod problem (see Proposition B.1 in [1]) that, for all $t \geq 0$,

$$(4.20) \quad \hat{I}_i^r(t) \geq -\inf_{0 \leq s \leq t}(M_{i1}^r\hat{X}_1^r(s) + M_{i2}^r\hat{X}_2^r(s) + M_{i3}^r\hat{X}_3^r(s)), \qquad i = 1, 2.$$

For $x \in \mathcal{D}^1$ with $x(0) = 0$, define $\Gamma(x) \in \mathcal{D}^1$ as

$$(4.21) \qquad\qquad \Gamma(x)(t) \doteq x(t) - \inf_{0 \leq s \leq t} x(s), \qquad t \geq 0.$$

Then (4.20) implies that

$$
\begin{aligned}
(4.22) \quad \hat{W}_i^r(t) &\geq \Gamma(M_{i1}^r\hat{X}_1^r(\cdot) + M_{i2}^r\hat{X}_2^r(\cdot) + M_{i3}^r\hat{X}_3^r(\cdot))(t) \\
&\qquad\qquad\qquad\qquad\qquad \text{for all } t \geq 0, i = 1, 2.
\end{aligned}
$$

Now we prove the inequality in (4.12). If the left-hand side of (4.12) is infinite, then the inequality holds trivially. Otherwise, get a further subsequence indexed by $r''$ (which may depend on $\omega, t$) such that

$$(4.23) \qquad \lim_{r'' \to \infty} h \cdot \hat{Q}^{r''}(t, \omega) = \liminf_{r' \to \infty} h \cdot \hat{Q}^{r'}(t, \omega) < \infty.$$

Notice that from (4.23), since $h_i > 0, \hat{Q}_i^{r''}(t, \omega) \geq 0$, for $i = 1, 2, 3$, we have $\{\hat{Q}_i^{r''}(t, \omega)\}$ is a bounded sequence as $r'' \to \infty$ for $i = 1, 2, 3$. From (4.14) and (4.17), we have

$$(4.24) \qquad \lim_{r'' \to \infty} (h_i^{r'', k(\omega, t)} - h_i) \hat{Q}_i^{r''}(t, \omega) = 0 \qquad \text{for } i = 1, 2, 3.$$

Using the above equality along with (4.19), (4.22) and the nonnegativity of $a_i^*, b_i^*; i = 1, 2$, we have

$$
\begin{aligned}
\lim_{r'' \to \infty} &(h \cdot \hat{Q}^{r''}(t, \omega)) \\
&= \lim_{r'' \to \infty} (h^{r'', k(\omega, t)} \cdot \hat{Q}^{r''}(t, \omega)) \\
&\geq \limsup_{r'' \to \infty} (a_{k(\omega, t)}^* \hat{W}_1^{r''}(t, \omega) + b_{k(\omega, t)}^* \hat{W}_2^{r''}(t, \omega)) \\
(4.25) \qquad &\geq \limsup_{r'' \to \infty} [a_{k(\omega, t)}^* \Gamma(M_{11}^{r''} \hat{X}_1^{r''}(\cdot, \omega) + M_{12}^{r''} \hat{X}_2^{r''}(\cdot, \omega) + M_{13}^{r''} \hat{X}_3^{r''}(\cdot, \omega))(t) \\
&\qquad\qquad + b_{k(\omega, t)}^* \Gamma(M_{21}^{r''} \hat{X}_1^{r''}(\cdot, \omega) + M_{22}^{r''} \hat{X}_2^{r''}(\cdot, \omega) + M_{23}^{r''} \hat{X}_3^{r''}(\cdot, \omega))(t)] \\
&= a_{k(\omega, t)}^* \Gamma(M_{11} \tilde{X}_1(\cdot, \omega) + M_{12} \tilde{X}_2(\cdot, \omega) + M_{13} \tilde{X}_3(\cdot, \omega))(t) \\
&\qquad + b_{k(\omega, t)}^* \Gamma(M_{21} \tilde{X}_1(\cdot, \omega) + M_{22} \tilde{X}_2(\cdot, \omega) + M_{23} \tilde{X}_3(\cdot, \omega))(t) \\
&= a_{k(\omega, t)}^* \tilde{W}_1^*(t, \omega) + b_{k(\omega, t)}^* \tilde{W}_2^*(t, \omega) \\
&= h \cdot \hat{Q}^*(t, \omega),
\end{aligned}
$$

where the fifth line follows from (4.10) and the continuity of the one-dimensional Skorohod map $\Gamma(\cdot)$ on $\mathcal{D}^1$ and recalling that $M^r \to M$ as $r \to \infty$. The last equality in (4.25) follows from the definition of $a_i^*, b_i^*, i = 1, 2$, definition of $k(\cdot, \cdot)$ and (3.23)–(3.24). This proves (4.12) and the result follows.  $\square$

We now proceed to the proof of Theorem 4.2. We begin with the following elementary result.

Lemma 4.7.  Let $\{f_r\}, \{g_r\}$ be sequences of functions in $\mathcal{D}^1$, and $f, g$ be functions in $\mathcal{C}^1$ such that $f_r \to f, g_r \to g$ in $\mathcal{D}^1$ as $r \to \infty$. Suppose that

$$(4.26) \qquad \int_0^\infty e^{-\gamma t} \mathbf{I}_{\{|g(t)|=0\}} \, dt = 0.$$

Let $\{\varepsilon_r\}$ be a sequence of nonnegative numbers such that $\varepsilon_r \to 0$ as $r \to \infty$. Then, for all $T > 0$, the following hold:

$$(4.27) \quad \int_0^T e^{-\gamma t} f_r(t) \mathbf{I}_{\{g_r(t) \geq \varepsilon_r\}} \, dt \to \int_0^T e^{-\gamma t} f(t) \mathbf{I}_{\{g(t) \geq 0\}} \, dt \qquad \text{as } r \to \infty,$$

$$(4.28) \quad \int_0^T e^{-\gamma t} f_r(t) \mathbf{I}_{\{g_r(t) < \varepsilon_r\}} \, dt \to \int_0^T e^{-\gamma t} f(t) \mathbf{I}_{\{g(t) \leq 0\}} \, dt \qquad \text{as } r \to \infty.$$

PROOF. Let $\mu$ be a finite measure defined on $(\mathbb{R}_+, \mathcal{B}(\mathbb{R}_+))$ via the following relation:

$$(4.29) \qquad\qquad d\mu(t) = e^{-\gamma t} \, dt.$$

It follows from (4.26) that

$$(4.30) \qquad\qquad g(t) \neq 0 \qquad \text{a.e. } t[\mu].$$

To prove (4.27), we will show that

$$(4.31) \quad \left| \int_0^T f_r(t) \mathbf{I}_{\{g_r(t) \geq \varepsilon_r\}} \, d\mu(t) - \int_0^T f(t) \mathbf{I}_{\{g(t) \geq 0\}} \, d\mu(t) \right| \to 0 \qquad \text{as } r \to \infty.$$

We can bound the left-hand side of (4.31) by

$$(4.32) \quad \int_0^T |f_r(t) - f(t)| \, d\mu(t) + \int_0^T |f(t)| \cdot |\mathbf{I}_{\{g_r(t) \geq \varepsilon_r\}} - \mathbf{I}_{\{g(t) \geq 0\}}| \, d\mu(t).$$

Since $f$ is continuous and $f_r \to f$ in $\mathcal{D}^1$, we have $\sup_{0 \leq t \leq T} |f_r(t) - f(t)| \to 0$. This shows that the first term in (4.32) converges to zero as $r \to \infty$.

For the second term in (4.32), it is enough to show the following:

$$(4.33) \qquad \mathbf{I}_{\{g_r(t) \geq \varepsilon_r\}} \to \mathbf{I}_{\{g(t) \geq 0\}} \qquad \text{as } r \to \infty \text{ for a.e. } t[\mu].$$

But (4.33) is an immediate consequence of (4.30) and the fact that, since $g$ is continuous,

$$(4.34) \quad \mathbf{I}_{\{g_r(t) \geq \varepsilon_r\}} \to \mathbf{I}_{\{g(t) \geq 0\}} \qquad \text{as } r \to \infty \text{ for all } t \text{ such that } g(t) \neq 0.$$

This proves (4.33) and completes the proof of (4.27). The proof of (4.28) is similar. $\qquad \square$

The following three theorems, the proofs of which are deferred to Section 5, are key to the proof of Theorem 4.2.

Let $\{T^r\}$ be the sequence of scheduling controls described in Definition 3.6. Let $\kappa$ be a positive constant satisfying

$$(4.35) \qquad \kappa > \max\left\{ 2\frac{\mu_1}{\mu_2}, 4, \frac{c}{(c-1)}, \frac{2\mu_2 c}{\mu_1(c-1)}, \theta_3 \right\},$$

where $\theta_3$ is as in Remark 4.3. For $r \in \mathbb{S}, t \geq 0$, define an event $\mathcal{E}(r,t)$ as follows:

$$(4.36) \quad \begin{aligned} \mathcal{E}(r,t) \doteq & \left\{ \sup_{0 \leq s \leq t} \hat{Q}_3^r(s) \mathbf{I}_{\{\hat{Q}_3^r(s) - (\mu_2^r/\mu_1^r)\hat{Q}_1^r(s) < \mathbf{L}^r/r\}} > \frac{\kappa(\mathbf{C}^r - \mathbf{L}^r + 1)}{r} \right\} \\ & \cup \left\{ \sup_{0 \leq s \leq t} \hat{Q}_1^r(s) \mathbf{I}_{\{\hat{Q}_3^r(s) - (\mu_2^r/\mu_1^r)\hat{Q}_1^r(s) \geq \mathbf{L}^r/r\}} > \frac{\kappa(\mathbf{C}^r - \mathbf{L}^r + 1)}{r} \right\}. \end{aligned}$$

Note that the event $\mathcal{E}(r, t)$ depends on parameters $c$ and $\ell_0$, however, this dependence is suppressed in the notation.

THEOREM 4.8.  *Let $\{T^r\}$ be the sequence of scheduling controls described in Definition 3.6 with threshold parameters $c$ and $\ell_0$. Let $\theta_3$ be as in Remark 4.3(a). Then there exist constants $\theta_i \in (0, \infty)$, $i = 1, 2$, and $r_0 \geq 1$ such that for all $t \in [0, \infty)$, $r \geq r_0$; $\ell_0 \in (1, \infty)$ and $c \geq 1 + \frac{4}{\theta_3}$,*

$$(4.37) \qquad \mathbb{P}(\mathcal{E}(r, t)) \leq \theta_1 (1 + r^4 t^2)(e^{-\theta_2 r^2 t} + r^{-\theta_3 (c-1)\ell_0}).$$

THEOREM 4.9.  *There exists $c \in [1 + \frac{4}{\theta_3}, \infty)$ such that, if $\{T^r\}$ is the sequence of scheduling controls described in Definition 3.6 with threshold parameters $c$ and some $\ell_0 \in (1, \infty)$, then there exist constants $\gamma_i > 0, i = 1, \ldots, 4$, $r_1 \geq 1$, $d \in (0, \infty)$ such that, for all $r \geq r_1$, $t \in [0, \infty)$, we have*

$$(4.38) \qquad \begin{aligned} &\mathbb{P}\left[\int_{[0,t)} \mathbf{I}_{\{\hat{Q}_2^r(s) \geq d\ell_0 \log r / r\}} \, d\hat{I}_2^r(s) \neq 0\right] \\ &\qquad \leq \gamma_1 (1 + r^2 t)e^{-\gamma_2 r^2 t} + \gamma_3 (1 + r^2 t)^2 r^{-\gamma_4 \ell_0}. \end{aligned}$$

Proofs of Theorems 4.8 and 4.9 will be given in Section 5. An immediate corollary of the above theorems is the following.

COROLLARY 4.10.  *Let $c$ and the scheduling sequence $\{T^r\}$ be as in Theorem 4.9. Suppose that $\bar{\ell} \in (0, \infty)$ is large enough so that $\theta_3 (c-1)\bar{\ell} > 4$ and $\gamma_4 \bar{\ell} > 4$. Then, for each fixed $t \geq 0$, for all $\ell_0 > \bar{\ell}$, the probabilities (4.37) and (4.38) tend to zero as $r \to \infty$. This, in particular, implies that, as $r \to \infty$,*

$$(4.39) \qquad \begin{aligned} &\hat{Q}_1^r(\cdot)\mathbf{I}_{\{\hat{Q}_3^r(\cdot) - (\mu_2^r / \mu_1^r)\hat{Q}_1^r(\cdot) \geq \mathbf{L}^r / r\}} \Rightarrow 0, \\ &\hat{Q}_3^r(\cdot)\mathbf{I}_{\{\hat{Q}_3^r(\cdot) - (\mu_2^r / \mu_1^r)\hat{Q}_1^r(\cdot) < \mathbf{L}^r / r\}} \Rightarrow 0, \\ &\int_{[0,\cdot]} \mathbf{I}_{\{\hat{Q}_2^r(s) \geq d\ell_0 \log r / r\}} \, d\hat{I}_2^r(s) \Rightarrow 0. \end{aligned}$$

Using the third convergence result above, we will obtain in Theorem 4.2(a) that $(\hat{W}^r, \hat{I}^r) \Rightarrow (\tilde{W}^*, \tilde{I}^*)$ as $r \to \infty$. However, we are unable to show that $\hat{Q}^r \Rightarrow \tilde{Q}^*$ as $r \to \infty$. Nevertheless, as will be seen in the proof of Theorem 4.8 below, the weak convergence results in Corollary 4.10 with suitable uniform integrability estimates (Theorem 4.11, see Remark 4.12) will suffice for the proof of asymptotic optimality of the proposed policy.

THEOREM 4.11.  *Suppose that $c$ is as obtained through Theorem 4.9 and $\ell_0$ satisfies the conditions in Corollary 4.10. Let $\{T^r\}$ be the sequence of* *scheduling controls described in Definition* 3.6 *with threshold parameters c and $\ell_0$. Then the following hold:*

$$(4.40) \qquad \limsup_{r \to \infty} \int_0^\infty e^{-\gamma t} \mathbb{E}\left[\sup_{0 \le s \le t} \hat{W}_i^r(s)\right]^2 dt < \infty, \qquad i = 1, 2.$$

*Also,*

$$(4.41) \quad \limsup_{T \to \infty} \limsup_{r \to \infty} \int_T^\infty e^{-\gamma t} \mathbb{E}\left[\sup_{0 \le s \le t} \hat{W}_i^r(s)\right]^2 dt = 0, \qquad i = 1, 2.$$

Proof of Theorem 4.11 will be given in Section 5.

REMARK 4.12. Note that (4.40) implies that $\hat{W}_i^r(\cdot)$, $i = 1, 2$, are uniformly integrable (u.i.) with respect to the product measure $\mathbb{P} \times \mu$, where $\mu$ is as defined in (4.29). Also, note that (4.41), in particular, implies that

$$(4.42) \qquad \limsup_{T \to \infty} \limsup_{r \to \infty} \int_T^\infty e^{-\gamma t} \mathbb{E}(\hat{W}_i^r(t)) \, dt = 0, \qquad i = 1, 2.$$

Lemma 4.6 showed that if $\{T^r\}$ is any sequence of admissible controls which gives a finite cost asymptotically, that is, (4.5) is satisfied, then $\bar{T}^r \Rightarrow \bar{T}^*$. For the sequence of scheduling controls in Definition 3.6, we do not know a priori that (4.5) is satisfied. In view of that, we prove the following lemma.

LEMMA 4.13. *Let $\{T^r\}$ be the sequence of scheduling controls described in Definition* 3.6 *with threshold parameters c and $\ell_0$. Suppose that c is as chosen in Theorem* 4.9 *and $\ell_0$ satisfies the conditions in Corollary* 4.10. *Then*

$$(4.43) \qquad\qquad\qquad \bar{T}^r \Rightarrow \bar{T}^*,$$

*where $\bar{T}^*$ is as defined in* (3.1).

PROOF. From (2.16), we have that, for $s \ge 0$,

$$(4.44) \qquad\qquad \hat{W}_1^r(s) = \frac{\hat{X}_1^r(s)}{\mu_1^r} + \frac{\hat{X}_2^r(s)}{\mu_2^r} + \hat{I}_1^r(s).$$

Now, the left-hand side of (4.44) is nonnegative, and by definition of the proposed scheduling policy (Definition 3.6), $\hat{I}_1^r$ is nondecreasing, starts from zero and increases only when both $\hat{Q}_1^r$ and $\hat{Q}_2^r$ are zero, or, in other words, in view of (2.10), only when the left-hand side of (4.44) is zero. Let $\Gamma(\cdot)$ be the Skorohod map defined in (4.21). From a well-known characterization of the solution of a one-dimensional Skorohod problem (see [1], Proposition B.1), we have from (4.44) that (4.20) and (4.22) hold with inequalities replaced by equalities. In particular,

$$\hat{W}_1^r(s) = \Gamma\left(\frac{\hat{X}_1^r(\cdot)}{\mu_1^r} + \frac{\hat{X}_2^r(\cdot)}{\mu_2^r}\right)(s) \tag{4.45}$$

and

$$\hat{I}_1^r(s) = -\inf_{0 \le u \le s}\left(\frac{\hat{X}_1^r(u)}{\mu_1^r} + \frac{\hat{X}_2^r(u)}{\mu_2^r}\right). \tag{4.46}$$

Next recall that

$$\hat{W}_2^r(s) = \frac{\hat{Q}_2^r(s)}{\mu_3^r} + \frac{\hat{Q}_3^r(s)}{\mu_3^r} = \frac{\hat{X}_2^r(s)}{\mu_3^r} + \frac{\hat{X}_3^r(s)}{\mu_3^r} + \hat{I}_2^r(s), \qquad s \ge 0. \tag{4.47}$$

From the second equality in (4.47), we get that, for $s \ge 0$,

$$
\begin{aligned}
\frac{\hat{Q}_2^r(s)}{\mu_3^r}&\mathbf{I}_{\{\hat{Q}_2^r(s) \ge 2d_0 \log r/r\}} + \frac{\hat{Q}_3^r(s)}{\mu_3^r} \\
&= \Big(\frac{\hat{X}_2^r(s)}{\mu_3^r} + \frac{\hat{X}_3^r(s)}{\mu_3^r} - \frac{\hat{Q}_2^r(s)}{\mu_3^r}\mathbf{I}_{\{\hat{Q}_2^r(s) < 2d_0 \log r/r\}} \\
&\qquad + \int_{[0,s]} \mathbf{I}_{\{\hat{Q}_2^r(u) \ge d_0 \log r/r\}}\, d\hat{I}_2^r(u)\Big) \\
&\qquad + \int_{[0,s]} \mathbf{I}_{\{\hat{Q}_2^r(u) < d_0 \log r/r\}}\, d\hat{I}_2^r(u).
\end{aligned}
\tag{4.48}
$$

Once more, from the definition of the scheduling policy (in Definition 3.6), the left-hand side of the equation in (4.48) is nonnegative and the last term on the right-hand side of (4.48) is nondecreasing, starts from zero and increases only when the left-hand side is zero. Also, note that since the paths of $\hat{Q}_2^r(\cdot)$ are piecewise constant, the processes $\frac{\hat{Q}_2^r(\cdot)}{\mu_3^r}\mathbf{I}_{\{\hat{Q}_2^r(\cdot) \ge 2d_0 \log r/r\}}$ and $\frac{\hat{Q}_2^r(\cdot)}{\mu_3^r}\mathbf{I}_{\{\hat{Q}_2^r(\cdot) < 2d_0 \log r/r\}}$ have paths in $\mathcal{D}^1$. Thus, using the characterizing property of the one-dimensional Skorohod map and (4.47)–(4.48), we obtain

$$
\begin{aligned}
\hat{W}_2^r(t) = \Gamma\Big(&\frac{\hat{X}_2^r(\cdot)}{\mu_3^r} + \frac{\hat{X}_3^r(\cdot)}{\mu_3^r} - \frac{\hat{Q}_2^r(\cdot)}{\mu_3^r}\mathbf{I}_{\{\hat{Q}_2^r(\cdot) < 2d_0 \log r/r\}} \\
&\qquad + \int_{[0,\cdot]} \mathbf{I}_{\{\hat{Q}_2^r(u) \ge d_0 \log r/r\}}\, d\hat{I}_2^r(u)\Big)(t) \\
&+ \frac{\hat{Q}_2^r(t)}{\mu_3^r}\mathbf{I}_{\{\hat{Q}_2^r(t) < 2d_0 \log r/r\}},
\end{aligned}
\tag{4.49}
$$

$$\int_{[0,s]} \mathbf{I}_{\{\hat{Q}_2^r(u) < d_0 \log r / r\}} \, d\hat{I}_2^r(u)$$

$$(4.50) \qquad = -\inf_{0 \leq s \leq t} \left( \frac{\hat{X}_2^r(s)}{\mu_3^r} + \frac{\hat{X}_3^r(s)}{\mu_3^r} - \frac{\hat{Q}_2^r(s)}{\mu_3^r} \mathbf{I}_{\{\hat{Q}_2^r(s) < 2d_0 \log r / r\}} \right.$$
$$\left. + \int_{[0,s]} \mathbf{I}_{\{\hat{Q}_2^r(u) \geq d_0 \log r / r\}} \, d\hat{I}_2^r(u) \right).$$

Using the fact that $\Gamma(\cdot)$ is Lipschitz continuous with constant 2 along with (4.45) and (4.49), we have the following:

$$(4.51) \qquad \begin{aligned} \sup_{0 \leq s \leq t} \bar{W}_1^r(s) &\leq \sup_{0 \leq s \leq t} \frac{1}{r} \Gamma\left( \frac{\hat{X}_1^r(\cdot)}{\mu_1^r} + \frac{\hat{X}_2^r(\cdot)}{\mu_2^r} \right)(s) \\ &\leq 2\frac{1}{r\mu_1^r} \sup_{0 \leq s \leq t} |\hat{X}_1^r(s)| + 2\frac{1}{r\mu_2^r} \sup_{0 \leq s \leq t} |\hat{X}_2^r(s)| \end{aligned}$$

and

$$\sup_{0 \leq s \leq t} \bar{W}_2^r(s)$$

$$\leq \sup_{0 \leq s \leq t} \frac{1}{r} \Gamma\left( \frac{\hat{X}_2^r(\cdot)}{\mu_3^r} + \frac{\hat{X}_3^r(\cdot)}{\mu_3^r} - \frac{\hat{Q}_2^r(\cdot)}{\mu_3^r} \mathbf{I}_{\{\hat{Q}_2^r(\cdot) < 2d_0 \log r / r\}} \right.$$
$$\left. + \int_{[0,\cdot]} \mathbf{I}_{\{\hat{Q}_2^r(u) \geq d_0 \log r / r\}} \, d\hat{I}_2^r(u) \right)(s)$$

$$(4.52) \qquad + \sup_{0 \leq s \leq t} \frac{1}{r} \frac{\hat{Q}_2^r(s)}{\mu_3^r} \mathbf{I}_{\{\hat{Q}_2^r(s) < 2d_0 \log r / r\}}$$

$$\leq 2\frac{1}{r\mu_2^r} \sup_{0 \leq s \leq t} |\hat{X}_2^r(s)| + 2\frac{1}{r\mu_3^r} \sup_{0 \leq s \leq t} |\hat{X}_3^r(s)|$$

$$+ 3\frac{1}{r} \sup_{0 \leq s \leq t} \frac{\hat{Q}_2^r(s)}{\mu_3^r} \mathbf{I}_{\{\hat{Q}_2^r(s) < 2d_0 \log r / r\}}$$

$$+ 2\frac{1}{r} \int_{[0,t]} \mathbf{I}_{\{\hat{Q}_2^r(u) \geq d_0 \log r / r\}} \, d\hat{I}_2^r(u).$$

From (3.3), (2.12) and Assumption 2.1, it follows that, for $i = 1, 2, 3$,

$$(4.53) \qquad \frac{1}{r\mu_i^r} \sup_{0 \leq s \leq t} |\hat{X}_i^r(s)| \to 0 \qquad \text{in probability, as } r \to \infty.$$

Also, note that

$$(4.54) \qquad \sup_{0 \leq s \leq t} \frac{\hat{Q}_2^r(s)}{\mu_3^r} \mathbf{I}_{\{\hat{Q}_2^r(s) < 2d_0 \log r / r\}} \leq \frac{2d_0 \log r}{r\mu_3^r} \to 0.$$

Now using (4.53), (4.54) and Corollary 4.10 in (4.51) and (4.52), we get that

$$(4.55) \qquad \bar{W}_i^r \Rightarrow 0, \qquad i = 1, 2.$$

This immediately yields that

$$(4.56) \qquad \bar{Q}_i^r \Rightarrow 0, \qquad i = 1, 2, 3.$$

Also, from (3.3), we have that

$$(4.57) \qquad \bar{A}_i^r(\cdot) \Rightarrow 0, \qquad \bar{S}_j^r(\bar{T}_j^r(\cdot)) \Rightarrow 0, \qquad i = 1, 2, j = 1, 2, 3.$$

Using (4.56) and (4.57), it follows from (2.11) that, for all $t \geq 0$,

$$(4.58) \qquad \begin{array}{l} \lambda_i^r t - \mu_i^r \bar{T}_i^r(t) \Rightarrow 0, \qquad i = 1, 2, \\ \mu_2^r \bar{T}_2^r(t) - \mu_3^r \bar{T}_3^r(t) \Rightarrow 0. \end{array}$$

The result follows on combining (4.58) with Assumption 2.1 and (2.14) of Assumption 2.2. $\square$

We now come to the proof of the main result of this section.

PROOF OF THEOREM 4.2. Suppose that $c$ is obtained from Theorem 4.9 and $\bar{\ell} \doteq \max\{\frac{4}{\theta_3(c-1)}, \frac{4}{\gamma_4}\}$, where $\theta_i, \gamma_i$ are as in Theorems 4.8 and 4.9, respectively. Henceforth, the sequence $\{T^r\}$ will have threshold parameters $c$ and $\ell_0$, with $\ell_0 \in (\bar{\ell}, \infty)$.

From Lemma 4.13 and (3.3), we have that

$$(4.59) \qquad \hat{X}_i^r \Rightarrow \tilde{X}_i \qquad \text{as } r \to \infty, i = 1, 2, 3,$$

where $\tilde{X}$ is as defined below (3.4). From this, Corollary 4.10, Assumption 2.1 and alternative expressions for $\hat{W}_i^r$, $\hat{I}_i^r$, for $i = 1, 2$, in (4.45), (4.46), (4.49) and (4.50), it follows that

$$(4.60) \qquad (\hat{W}^r, \hat{I}^r) \Rightarrow (\tilde{W}^*, \tilde{I}^*) \qquad \text{as } r \to \infty.$$

We have also used (4.54) and continuity of $\Gamma(\cdot)$ in obtaining (4.60). This proves part (a) of the theorem.

For part (b), first we observe that from Theorem 4.11 (see Remark 4.12) and part (a) of this theorem, it follows that

$$(4.61) \qquad \int_0^\infty e^{-\gamma t} \mathbb{E}(\hat{W}_i^r(t)) \, dt \to \int_0^\infty e^{-\gamma t} \mathbb{E}(\tilde{W}_i^*(t)) \, dt, \qquad i = 1, 2.$$

Next observe that the reflected Brownian motions $\tilde{W}_1^*$ and $\tilde{W}_2^*$ satisfy, for every $t \geq 0$, $\mathbb{P}(\mu_3 \tilde{W}_2^*(t) = \mu_2 \tilde{W}_1^*(t)) = 0$. Using this fact and Fubini's theorem, it follows that

$$(4.62) \qquad \int_0^\infty e^{-\gamma t} \mathbf{I}_{\{\mu_3 \tilde{W}_2^*(t) - \mu_2 \tilde{W}_1^*(t) = 0\}} \, dt = 0 \qquad \text{a.s. } [\mathbb{P}].$$

From Lemma 4.7 [see (4.27)], (4.60), (4.62) and the fact that $\frac{\mathbf{L}^r}{r} = \ell_0 \frac{\log r}{r}$ decreases to 0 as $r \to \infty$, we have, for all $T \geq 0$,

$$
(4.63) \qquad \begin{aligned}
\int_0^T e^{-\gamma t} \hat{W}_1^r(t) \mathbf{I}_{\{\mu_3^r \hat{W}_2^r(t) - \mu_2^r \hat{W}_1^r(t) \geq \mathbf{L}^r/r\}} \, dt \\
\to \int_0^T e^{-\gamma t} \tilde{W}_1^*(t) \mathbf{I}_{\{\mu_3 \tilde{W}_2^*(t) - \mu_2 \tilde{W}_1^*(t) \geq 0\}} \, dt
\end{aligned}
$$

in distribution. Using uniform integrability of $\hat{W}_1^r$ (see Remark 4.12), we can conclude from (4.63) that, for all $T \geq 0$,

$$
(4.64) \qquad \begin{aligned}
\int_0^T e^{-\gamma t} \mathbb{E}(\hat{W}_1^r(t) \mathbf{I}_{\{\mu_3^r \hat{W}_2^r(t) - \mu_2^r \hat{W}_1^r(t) \geq \mathbf{L}^r/r\}}) \, dt \\
\to \int_0^T e^{-\gamma t} \mathbb{E}(\tilde{W}_1^*(t) \mathbf{I}_{\{\mu_3 \tilde{W}_2^*(t) - \mu_2 \tilde{W}_1^*(t) \geq 0\}}) \, dt.
\end{aligned}
$$

From (4.42), for $i = 1$, and (4.64), simple calculations show that

$$
(4.65) \qquad \begin{aligned}
\int_0^\infty e^{-\gamma t} \mathbb{E}(\hat{W}_1^r(t) \mathbf{I}_{\{\mu_3^r \hat{W}_2^r(t) - \mu_2^r \hat{W}_1^r(t) \geq \mathbf{L}^r/r\}}) \, dt \\
\to \int_0^\infty e^{-\gamma t} \mathbb{E}(\tilde{W}_1^*(t) \mathbf{I}_{\{\mu_3 \tilde{W}_2^*(t) - \mu_2 \tilde{W}_1^*(t) \geq 0\}}) \, dt.
\end{aligned}
$$

Similarly, using (4.28) of Lemma 4.7 and (4.42) for $i = 2$, it can be shown that

$$
(4.66) \qquad \begin{aligned}
\int_0^\infty e^{-\gamma t} \mathbb{E}(\hat{W}_2^r(t) \mathbf{I}_{\{\mu_3^r \hat{W}_2^r(t) - \mu_2^r \hat{W}_1^r(t) < \mathbf{L}^r/r\}}) \, dt \\
\to \int_0^\infty e^{-\gamma t} \mathbb{E}(\tilde{W}_2^*(t) \mathbf{I}_{\{\mu_3 \tilde{W}_2^*(t) - \mu_2 \tilde{W}_1^*(t) \leq 0\}}) \, dt.
\end{aligned}
$$

From (4.65), (4.66), (3.23)–(3.24) and Assumption 2.1, it follows that

$$
(4.67) \qquad \begin{aligned}
\int_0^\infty e^{-\gamma t} \mathbb{E}(\mu_2^r \hat{W}_1^r(t) \mathbf{I}_{\{\mu_3^r \hat{W}_2^r(t) - \mu_2^r \hat{W}_1^r(t) \geq \mathbf{L}^r/r\}}) \, dt \\
+ \int_0^\infty e^{-\gamma t} \mathbb{E}(\mu_3^r \hat{W}_2^r(t) \mathbf{I}_{\{\mu_3^r \hat{W}_2^r(t) - \mu_2^r \hat{W}_1^r(t) < \mathbf{L}^r/r\}}) \, dt \\
\to \int_0^\infty e^{-\gamma t} \mathbb{E}(\tilde{Q}_2^*(t)) \, dt.
\end{aligned}
$$

Now using (2.10), the left-hand side of (4.67) can be written as

$$
(4.68) \qquad \begin{aligned}
\int_0^\infty e^{-\gamma t} \mathbb{E}(\hat{Q}_2^r(t)) \, dt \\
+ \int_0^\infty e^{-\gamma t} \mathbb{E}\left( \frac{\mu_2^r \hat{Q}_3^r(s)}{\mu_1^r} \mathbf{I}_{\{\hat{Q}_3^r(s) - (\mu_2^r/\mu_1^r)\hat{Q}_1^r(s) \geq \mathbf{L}^r/r\}} \right) dt \\
+ \int_0^\infty e^{-\gamma t} \mathbb{E}(\hat{Q}_3^r(s) \mathbf{I}_{\{\hat{Q}_3^r(s) - (\mu_2^r/\mu_1^r)\hat{Q}_1^r(s) < \mathbf{L}^r/r\}}) \, dt.
\end{aligned}
$$

From the uniform integrability of $\hat{W}_i^r$ given in Remark 4.12 and recalling that $\mu_i^r \to \mu_i, i = 1, 2$, we have that, for $j = 1, 2, 3$, $\hat{Q}_j^r$ are uniformly integrable (with respect to the measure $\mathbb{P} \times \mu$). Combining this observation with Corollary 4.10, it follows that the last two terms of (4.68) tend to zero. This, in view of (4.67) and Assumption 2.1, implies that

$$(4.69) \qquad \int_0^\infty e^{-\gamma t} \mathbb{E}(\hat{Q}_2^r(t)) \, dt \to \int_0^\infty e^{-\gamma t} \mathbb{E}(\tilde{Q}_2^*(t)) \, dt \qquad \text{as } r \to \infty.$$

Now using (4.69) and (4.61) in (2.10) and (3.8), it follows immediately that

$$(4.70) \qquad \begin{aligned} &\int_0^\infty e^{-\gamma t} \mathbb{E}(\hat{Q}_i^r(t)) \, dt \\ &\qquad \to \int_0^\infty e^{-\gamma t} \mathbb{E}(\tilde{Q}_i^*(t)) \, dt \qquad \text{as } r \to \infty, i = 1, 3. \end{aligned}$$

Finally, combining the above two displays with the definition of the cost function and the representation of $J^*$ in (3.28), it follows that

$$(4.71) \qquad \hat{J}^r(T^r) = \sum_{i=1}^3 \mathbb{E}\left( \int_0^\infty e^{-\gamma t} h_i \hat{Q}_i^r(t) \, dt \right) \to J^* \qquad \text{as } r \to \infty.$$

This completes the proof of the theorem. $\square$

**5. Proofs of Theorems 4.8, 4.9 and 4.11.** We begin with the following standard large deviations estimate for Poisson processes. This estimate will be used in many of the arguments in this section. For a proof we refer the reader to [9] or Theorem 5.3 of [15].

THEOREM 5.1 (Kurtz [9]). *Let $N_\lambda(\cdot)$ be a Poisson process with rate $\lambda > 0$. Then for all $0 < \underline{\lambda} < \overline{\lambda} < \infty$, there exists a $\tilde{C}_1 \in (0, \infty)$ and a function $\tilde{C}_2 : (0, \infty) \to (0, \infty)$ such that, for all $\alpha > 0, \varepsilon > 0$,*

$$(5.1) \qquad \sup_{\lambda \in [\underline{\lambda}, \overline{\lambda}]} \mathbb{P}\left( \sup_{0 \le t \le 1} \left| \frac{N_\lambda(\alpha t)}{\alpha} - \lambda t \right| \ge \varepsilon \right) \le \tilde{C}_1 e^{-\alpha \tilde{C}_2(\varepsilon)}.$$

An immediate corollary of the above result is the following.

COROLLARY 5.2. *Let $\{N^r(\cdot)\}_{r \in \mathbb{S}}$ be a sequence of Poisson processes with rates $\varrho^r$ such that $\varrho^r \to \varrho \in (0, \infty)$ as $r \to \infty$. Then there exists a $C_1 \in (0, \infty)$ and a function $C_2 : (0, \infty) \to (0, \infty)$ such that, for all $\varepsilon' > 0, \theta \in (0, 1), c^* > 0, r \in \mathbb{S}$, we have*

$$(5.2) \qquad \mathbb{P}\left( \sup_{\theta c^* \log r \le s \le c^* \log r} \left| \frac{N^r(s) - \varrho^r s}{s} \right| \ge \varepsilon' \right) \le C_1 e^{-c^* C_2(\varepsilon' \theta) \log r}.$$

The above corollary follows from some straightforward calculations on setting $\alpha = c^* \log r$ and $\varepsilon = \theta \varepsilon'$ in Theorem 5.1. Another important consequence of Theorem 5.1 is the following "terminal time" estimate.

COROLLARY 5.3. *Let $\{N^r(\cdot)\}_{r\in\mathbb{S}}$ be as in Corollary 5.2. Let $\varepsilon > 0$ be arbitrary. Then there exist $\varsigma_i \equiv \varsigma_i(\{\varrho^r\}, \varepsilon) \in (0, \infty)$, $i = 1, 2$, and $r_1 \equiv r_1(\{\varrho^r\}, \varepsilon) \in (0, \infty)$ such that, for all $r \geq r_1$,*

$$(5.3) \quad \mathbb{P}(N^r(t) \geq (\varrho^r + \varepsilon)t \ or \ N^r(t) \leq (\varrho^r - \varepsilon)t) \leq \varsigma_1 e^{-\varsigma_2 t} \qquad \forall t \in [0, \infty).$$

REMARK 5.4. The large deviation estimates in the above three results are used in the proofs of Theorems 4.8 and 4.9, which, in turn, are key to the proof of Theorem 4.2. Proof of Theorem 4.1 does not rely on any large deviation estimates and can be extended in a straightforward manner to the case of more general inter-arrival and service time distributions for $\{u_k^r(i), v_j^r(i) : k = 1, 2; j = 1, 2, 3; i = 1, 2, \ldots\}$ such that the corresponding renewal processes, $A_k^r(\cdot)$, $S_j^r(\cdot)$, satisfy a functional central limit result similar to (3.3). However, in order to extend Theorem 4.2 to a general renewal process setting, more stringent moment conditions on the above distributions are needed. Proof of Theorem 4.8 uses the one-dimensional large deviation estimate in Corollary 5.3. Corresponding results for a general renewal process are well known and indeed were used by the authors in [1] to prove the asymptotic optimality of their policy. Under precisely the assumptions of [1] on the underlying renewal process (Assumption 3.3 of that paper), one can extend the proof of Theorem 4.8 to a nonexponential setting. Note, however, that the proof will need to be modified to account for the non-Markovity by using multi-parameter filtrations and stopping times and using Lemma 7.6 of [1] in place of the strong Markov property. These modifications are fairly straightforward. Proof of Theorem 4.9 crucially relies on Corollary 5.2, which is a statement on the sample path large deviations of the underlying renewal process. We conjecture that using Theorem 3.1 of [14], one can extend the proof of Theorem 4.9 to a larger class of renewal processes which satisfy suitable exponential moment conditions.

Now, we proceed to the proof of Theorem 4.8. We begin by defining the following family of stopping times with respect to the filtration $\{\mathcal{F}_t^r\}_{t\geq 0}$, where $\mathcal{F}_t^r \doteq \sigma\{Q_j^r(s) : 0 \leq s \leq t, j = 1, 2, 3\}$. For $r \in \mathbb{S}$ and $n = 1, 2, \ldots$, define

$$
\begin{aligned}
(5.4) \quad \tau_0^r &\doteq 0, \\
\tau_{2n-1}^r &= \inf\left\{ t > \tau_{2n-2}^r \mid Q_3^r(t) - \frac{\mu_2^r}{\mu_1^r} Q_1^r(t) \geq \mathbf{L}^r \right\}, \\
\tau_{2n}^r &= \inf\left\{ t > \tau_{2n-1}^r \mid Q_3^r(t) - \frac{\mu_2^r}{\mu_1^r} Q_1^r(t) < \mathbf{L}^r \ \text{and} \ Q_3^r(t) < \mathbf{C}^r - 1 \right\}.
\end{aligned}
$$

From the form of the scheduling policy in Definition 3.6, it follows that $Q_3^r(\tau_{2n-2}^r) < \mathbf{C}^r - 1$. Thus, $Q_3^r(s)$ starts from below $\mathbf{C}^r - 1$ on $[\tau_{2n-2}^r, \tau_{2n-1}^r)$, and whenever the queue-length crosses $(\mathbf{C}^r - 1)$, *Server* 1 stops serving *Buffer* 2, causing $Q_3^r(s)$ to decrease monotonically. Thus, we have that

$$(5.5) \qquad Q_3^r(s) \leq \mathbf{C}^r \qquad \text{for all } s \in [\tau_{2n-2}^r, \tau_{2n-1}^r], n = 1, 2, \ldots.$$

PROOF OF THEOREM 4.8.   Recalling the definition of the diffusion-scaled processes $\hat{Q}_1^r(s), \hat{Q}_3^r(s)$, we can rewrite $\mathcal{E}(r,t)$ from (4.36) as

$$
\begin{aligned}
(5.6) \quad & \left\{ \sup_{0 \le s \le r^2 t} Q_3^r(s) \mathbf{I}_{\{Q_3^r(s) - (\mu_2^r/\mu_1^r)Q_1^r(s) < \mathbf{L^r}\}} > \kappa(\mathbf{C^r} - \mathbf{L^r} + 1) \right\} \\
& \cup \left\{ \sup_{0 \le s \le r^2 t} Q_1^r(s) \mathbf{I}_{\{Q_3^r(s) - (\mu_2^r/\mu_1^r)Q_1^r(s) \ge \mathbf{L^r}\}} > \kappa(\mathbf{C^r} - \mathbf{L^r} + 1) \right\}.
\end{aligned}
$$

Let

$$
(5.7) \qquad\qquad n^r \doteq [(\lambda_1^r + \lambda_2^r + 2)r^2 t] + 1.
$$

Note that every $\tau_{2k-1}^r$ $(k = 1, 2, \dots)$ corresponds to one up-crossing of $Q_3^r(s) - \frac{\mu_2^r}{\mu_1^r} Q_1^r(s)$ from below $\mathbf{L^r}$ to the threshold level $\mathbf{L^r}$ or above. Each such up-crossing either requires at least one service of a Class 1 job, which, in turn, implies at least 1 arrival of a Class 1 job, or it requires one service of a Class 2 job, which implies 1 arrival of a Class 2 job has occurred. Thus, the number of $\tau_{2k-1}^r$ in the interval $[0, r^2 t]$ is bounded above by $A_1^r(r^2 t) + A_2^r(r^2 t)$. Therefore,

$$
\begin{aligned}
\mathbb{P}(\tau_{2n^r-1}^r \le r^2 t) & \le \mathbb{P}(A_1^r(r^2 t) + A_2^r(r^2 t) \ge n^r) \\
& \le \mathbb{P}(A_1^r(r^2 t) + A_2^r(r^2 t) \ge (\lambda_1^r + \lambda_2^r + 2)r^2 t) \\
& \le \mathbb{P}(A_1^r(r^2 t) \ge (\lambda_1^r + 1)r^2 t) + \mathbb{P}(A_2^r(r^2 t) \ge (\lambda_2^r + 1)r^2 t) \\
& \le \kappa_1 e^{-\kappa_2 r^2 t},
\end{aligned}
$$

(5.8)

for all $r \ge \tilde{r}_1 \doteq \max\{r_1(\{\lambda_1^r\}, 1), r_1(\{\lambda_2^r\}, 1)\}$, where $\kappa_1 \doteq \varsigma_1(\{\lambda_1^r\}, 1) + \varsigma_1(\{\lambda_2^r\}, 1)$, $\kappa_2 \doteq \min\{\varsigma_2(\{\lambda_1^r\}, 1), \varsigma_2(\{\lambda_2^r\}, 1)\}$, and $r_1(\cdot)$, $\varsigma_1(\cdot)$, $\varsigma_2(\cdot)$ are as in Corollary 5.3.

Using (5.8) and the representation for $\mathcal{E}(r,t)$ in (5.6), we have that, for $r$ sufficiently large,

$$
\begin{aligned}
& \mathbb{P}(\mathcal{E}(r,t)) \\
& \quad \le \mathbb{P}(\tau_{2n^r-1}^r \le r^2 t) \\
& \qquad + \mathbb{P}\Big( \tau_{2n^r-1}^r > r^2 t, \\
& \qquad\qquad \left\{ \sup_{0 \le s \le r^2 t} Q_3^r(s) \mathbf{I}_{\{Q_3^r(s) - (\mu_2^r/\mu_1^r)Q_1^r(s) < \mathbf{L^r}\}} > \kappa(\mathbf{C^r} - \mathbf{L^r} + 1) \right\} \\
& \qquad\qquad \cup \left\{ \sup_{0 \le s \le r^2 t} Q_1^r(s) \mathbf{I}_{\{Q_3^r(s) - (\mu_2^r/\mu_1^r)Q_1^r(s) \ge \mathbf{L^r}\}} > \kappa(\mathbf{C^r} - \mathbf{L^r} + 1) \right\} \Big) \\
& \quad \le \kappa_1 e^{-\kappa_2 r^2 t} \\
& \qquad + \sum_{n=1}^{n^r} \mathbb{P}\Big( \tau_{2n-1}^r \le r^2 t, Q_3^r(s) > \kappa(\mathbf{C^r} - \mathbf{L^r} + 1)
\end{aligned}
$$

(5.9)

$$(5.10) \qquad \text{and } Q_3^r(s) - \frac{\mu_2^r}{\mu_1^r} Q_1^r(s) < \mathbf{L^r} \text{ for some } s \in [\tau_{2n-1}^r, \tau_{2n}^r \wedge (r^2 t)] \Big)$$

$$+ \sum_{n=1}^{n^r} \mathbb{P}\Big(\tau_{2n-1}^r \leq r^2 t, Q_1^r(s) > \kappa(\mathbf{C^r} - \mathbf{L^r} + 1)$$

$$(5.11) \qquad \text{and } Q_3^r(s) - (\mu_2^r/\mu_1^r) Q_1^r(s) \geq \mathbf{L^r}$$

$$\text{for some } s \in [\tau_{2n-1}^r, \tau_{2n}^r \wedge (r^2 t)] \Big)$$

$$\leq \kappa_1 e^{-\kappa_2 r^2 t}$$

$$(5.12) \qquad + 2 \sum_{n=1}^{n^r} \mathbb{P}(\tau_{2n-1}^r \leq r^2 t, Q_1^r(s) > \kappa''(\mathbf{C^r} - \mathbf{L^r} + 1)$$

$$\text{for some } s \in [\tau_{2n-1}^r, \tau_{2n}^r \wedge (r^2 t)]),$$

where $\kappa'' \doteq \min\{\kappa, \frac{\kappa \mu_1}{2\mu_2}\}$ and the display in (5.12) is a consequence of the fact that, for each summand in (5.10), $Q_3^r(s) - \frac{\mu_2^r}{\mu_1^r} Q_1^r(s) < \mathbf{L^r}$. Combining this with the condition $Q_3^r(s) > \kappa(\mathbf{C^r} - \mathbf{L^r} + 1)$ gives that $Q_1^r(s) > \frac{\mu_1^r}{\mu_2^r}[\kappa(\mathbf{C^r} - \mathbf{L^r} + 1) - \mathbf{L^r}] = \frac{3\kappa}{4}\frac{\mu_1^r}{\mu_2^r}(\mathbf{C^r} - \mathbf{L^r} + 1) + \frac{\mu_1^r}{\mu_2^r}[\frac{\kappa}{4}(\mathbf{C^r} - \mathbf{L^r} + 1) - \mathbf{L^r}] > \frac{3\kappa}{4}\frac{\mu_1^r}{\mu_2^r}(\mathbf{C^r} - \mathbf{L^r} + 1)$, using (4.35) and the fact that $c \geq 1 + 4/\theta_3$. Choosing $r$ to be sufficiently large, so that $\frac{\mu_1^r}{\mu_2^r} > \frac{2\mu_1}{3\mu_2}$ gives that, for such $r$, $Q_1^r(s) > \kappa''(\mathbf{C^r} - \mathbf{L^r} + 1)$. The sum (5.10) follows from (5.8), the fact that (4.35) implies $\mathbf{C^r} < \kappa(\mathbf{C^r} - \mathbf{L^r} + 1)$ and (5.5). The third term (5.11) is obtained using the fact that the indicator restricts us to the values of $s$ for which $Q_3^r(s) - \frac{\mu_2^r}{\mu_1^r} Q_1^r(s) \geq \mathbf{L^r}$, which happens only for $s \in [\tau_{2n-1}^r, \tau_{2n}^r)$.

Note that by our choice of $\kappa$ [see (4.35)], we have that

$$(5.13) \qquad \kappa'' \geq \max\left\{2\frac{\mu_1}{\mu_2}, \frac{c}{c-1}\right\}.$$

For fixed $r \in \mathbb{S}$ and $n \geq 1$, define a sequence of "intermediate" stopping times within $[\tau_{2n-1}^r, \tau_{2n}^r]$ as follows. For $m = 1, 2, \ldots,$

$$\eta_0^{r,n} \doteq \tau_{2n-1}^r,$$

$$\eta_{2m-1}^{r,n} \doteq \min\Big[\tau_{2n}^r, \inf\Big\{s > \eta_{2m-2}^{r,n} \Big| \Big(Q_3^r(s) - \frac{\mu_2^r}{\mu_1^r} Q_1^r(s) < \mathbf{L^r}$$

$$(5.14) \qquad\qquad \text{and } Q_3^r(s) \geq (\mathbf{C^r} - 1)\Big),$$

$$\text{OR } \Big(Q_3^r(s) - \frac{\mu_2^r}{\mu_1^r} Q_1^r(s) \geq \mathbf{L^r}$$

$$\text{and } Q_1^r(s) \geq \frac{\mu_1^r}{\mu_2^r}(\mathbf{C^r} - \mathbf{L^r} + 2)\Big)\Big\}\Big],$$

$$\eta_{2m}^{r,n} \doteq \min\Big[\tau_{2n}^r, \inf\Big\{s > \eta_{2m-1}^{r,n}\Big|Q_3^r(s) - \frac{\mu_2^r}{\mu_1^r}Q_1^r(s) \geq \mathbf{L^r}$$

$$\text{and } Q_1^r(s) < \frac{\mu_1^r}{\mu_2^r}(\mathbf{C^r} - \mathbf{L^r} + 2)\Big\}\Big].$$

Now, we estimate how many $\eta_{\cdot}^{r,n}$'s there can be in $[\tau_{2n-1}^r, \tau_{2n}^r \wedge (r^2 t)]$. Let $n^r$ be as in (5.7). Note that $\eta_{2n^r-1}^{r,n} \leq r^2 t$ implies that there are at least a total of $n^r$ arrivals in Class 1 and Class 2 together in $[0, r^2 t]$. Using an argument similar to that used in obtaining (5.8), we have that

$$(5.15) \qquad \mathbb{P}(\eta_{2n^r-1}^{r,n} \leq \tau_{2n}^r \wedge (r^2 t)) \leq \kappa_1 e^{-\kappa_2 r^2 t},$$

where $\kappa_i$ are as in (5.8). Now, each summand in (5.12) can be split over the sub-intervals formed by the $\eta_{\cdot}^{r,n}$'s spanning $[\tau_{2n-1}^r, \tau_{2n}^r \wedge (r^2 t)]$ as

$$\mathbb{P}(\tau_{2n-1}^r \leq r^2 t, Q_1^r(s) > \kappa''(\mathbf{C^r} - \mathbf{L^r} + 1) \text{ for some } s \in [\tau_{2n-1}^r, \tau_{2n}^r \wedge (r^2 t)])$$

$$(5.16) \qquad \leq \kappa_1 e^{-\kappa_2 r^2 t} + \sum_{m=1}^{n^r} \mathbb{P}\Big(\eta_{2m-1}^{r,n} \leq r^2 t, Q_1^r(\eta_{2m-1}^{r,n}) < \frac{\mu_1^r}{\mu_2^r}(\mathbf{C^r} - \mathbf{L^r} + 2) + 1,$$

$$Q_1^r(s) > \kappa''(\mathbf{C^r} - \mathbf{L^r} + 1)$$

$$\text{for some } s \in [\eta_{2m-1}^{r,n}, \eta_{2m}^{r,n} \wedge (r^2 t)]\Big).$$

Note that the terms corresponding to $s \in [\eta_{2m-2}^{r,n}, \eta_{2m-1}^{r,n} \wedge (r^2 t)]$ do not contribute to the sum (for large values of $r$), as the corresponding probabilities are zero. This is because, for $s \in [\eta_{2m-2}^{r,n}, \eta_{2m-1}^{r,n})$, $Q_1^r(s) < \frac{\mu_1^r}{\mu_2^r}(\mathbf{C^r} - \mathbf{L^r} + 2)$, which in view of (5.13), implies, for large $r$, $Q_1^r(s) < \kappa''(\mathbf{C^r} - \mathbf{L^r} + 1)$ for all $s \in [\eta_{2m-2}^{r,n}, \eta_{2m-1}^{r,n} \wedge (r^2 t)]$. Note that in all these calculations definition of $\tau_{2n}^r$ is used. This observation also provides the bound that if $\eta_{2m-1}^{r,n} \leq r^2 t$, we have $Q_1^r(\eta_{2m-1}^{r,n}) < \frac{\mu_1^r}{\mu_2^r}(\mathbf{C^r} - \mathbf{L^r} + 2) + 1$, which is used in (5.16). Combining (5.12) and (5.16), we have, for sufficiently large $r$ (choice of $r$ does not depend on $t$),

$$\mathbb{P}(\mathcal{E}(r,t))$$

$$\leq \kappa_1 e^{-\kappa_2 r^2 t}$$

$$(5.17) \qquad + 2\sum_{n=1}^{n^r}\Big[\kappa_1 e^{-\kappa_2 r^2 t} + \sum_{m=1}^{n^r} \mathbb{P}\Big(\eta_{2m-1}^{r,n} \leq r^2 t, Q_1^r(\eta_{2m-1}^{r,n})$$

$$\leq \frac{\mu_1^r}{\mu_2^r}(\mathbf{C^r} - \mathbf{L^r} + 2) + 1,$$

$$Q_1^r(s) > \kappa''(\mathbf{C^r} - \mathbf{L^r} + 1)$$

$$\text{for some } s \in [\eta_{2m-1}^{r,n}, \eta_{2m}^{r,n} \wedge (r^2 t)]\Big)\Big].$$

Let $\mathcal{A}_m^{r,n}$ denote the event in the $(m,n)$th summand in the last term of (5.17). On $\mathcal{A}_m^{r,n}$, our policy requires *Server* 1 to work continuously on *Buffer* 1 for all $s \in [\eta_{2m-1}^{r,n}, \eta_{2m}^{r,n} \wedge (r^2 t))$ and at the beginning of the interval, $Q_1^r$ is at most $\frac{\mu_1^r}{\mu_2^r}(\mathbf{C^r - L^r} + 2) + 1$. Also note that, for $r$ sufficiently large, *Buffer* 1 cannot be empty during this period. Indeed, for $s \in [\eta_{2m-1}^{r,n}, \eta_{2m}^{r,n} \wedge (r^2 t))$, we have $s < \eta_{2m}^{r,n} \leq \tau_{2n}^r$ and so by definition of $\tau_{2n}^r$, either $Q_3^r(t) - \frac{\mu_2^r}{\mu_1^r} Q_1^r(t) \geq \mathbf{L^r}$ or $Q_3^r(t) \geq \mathbf{C^r} - 1$ has to hold, and since $\eta_{2m}^{r,n} \geq \eta_{2m-1}^{r,n}$, by definition of $\eta_{2m}^{r,n}$, one of the following things must be true:

$$(5.18) \qquad Q_3^r(s) - \frac{\mu_2^r}{\mu_1^r} Q_1^r(s) < \mathbf{L^r} \quad \text{and} \quad Q_3^r(s) \geq \mathbf{C^r} - 1,$$

$$(5.19) \qquad Q_3^r(s) - \frac{\mu_2^r}{\mu_1^r} Q_1^r(s) \geq \mathbf{L^r} \quad \text{and} \quad Q_1^r(s) \geq \frac{\mu_1^r}{\mu_2^r}(\mathbf{C^r - L^r} + 2).$$

It is easy to see that in either case $Q_1^r(s) > 0$, using the fact that $\mathbf{C^r - L^r} \geq 1$ for $r$ sufficiently large.

Let, for $r \in \mathbb{S}$, $\mathcal{Q}^r(\cdot)$ denote an $M/M/1$ queue-length, with arrivals at rate $\lambda_1^r$ and service times at rate $\mu_1^r$. Define the stopping time

$$\beta^r \doteq \inf\left\{ s > 0 : \mathcal{Q}^r(s) < \frac{\mu_1^r}{\mu_2^r}(\mathbf{C^r - L^r} - 1) \right\}.$$

Then using the memoryless property of the exponential distribution and the form of the scheduling policy, it follows that each summand in (5.17) is bounded by

$$(5.20) \qquad \mathbb{P}\Big( \mathcal{Q}^r(s) > \kappa''(\mathbf{C^r - L^r} + 1) \text{ for some } s \in [0, \beta^r],$$
$$\mathcal{Q}^r(0) \leq \frac{\mu_1^r}{\mu_2^r}(\mathbf{C^r - L^r} + 2) + 1 \Big).$$

For $\varepsilon_2 > 0$, define $s^r \doteq [\frac{\mu_1^r}{\mu_2^r}(\mathbf{C^r - L^r} - 1) - 2]/2(\lambda_1^r + \varepsilon_2)$. Also define

$$(5.21) \qquad \mathcal{Y}^r \doteq \{\tilde{A}^r(s^r) < (\lambda_1^r + \varepsilon_2)s^r, \tilde{S}^r(s^r) > (\mu_1^r - \varepsilon_2)s^r \},$$

where $\tilde{A}^r$ and $\tilde{S}^r$ are the arrival and service processes of the $M/M/1$ queue. From Corollary 5.3, it follows that

$$(5.22) \qquad \mathbb{P}(\{\mathcal{Y}^r\}^c) \leq \kappa_3 e^{-\kappa_4 s^r}$$

for all $r \geq \tilde{r}_2 \doteq \max\{r_1(\{\lambda_1^r\}, \varepsilon_2), r_1(\{\mu_1^r\}, \varepsilon_2)\}$, where $\kappa_3 \doteq \varsigma_1(\{\lambda_1^r\}, \varepsilon_2) + \varsigma_1(\{\mu_1^r\}, \varepsilon_2)$, $\kappa_4 \doteq \min\{\varsigma_2(\{\lambda_1^r\}, \varepsilon_2), \varsigma_2(\{\mu_1^r\}, \varepsilon_2)\}$ and $r_1(\cdot)$, $\varsigma_1(\cdot)$, $\varsigma_2(\cdot)$ are as in Corollary 5.3.

Let $\tilde{\mathcal{A}}^r$ denote the event in (5.20). First, we argue that on $\mathcal{Y}^r \cap \tilde{\mathcal{A}}^r$, we have $s^r \geq \beta^r$. To show this, we argue by contradiction. Note that if $s^r < \beta^r$, then, since $\mathcal{Q}^r(s) \neq 0$ for all $s < \beta^r$, we have

$$\mathcal{Q}^r(s^r) = \mathcal{Q}^r(0) + \tilde{A}^r(s^r) - \tilde{S}^r(s^r)$$
$$\leq \frac{\mu_1^r}{\mu_2^r}(\mathbf{C^r} - \mathbf{L^r} + 2) + 1 + (\lambda_1^r - \mu_1^r + 2\varepsilon_2)s^r.$$

Choose $r_3$ large enough and $\varepsilon_2 > 0$ small enough, such that $(\lambda_1^r - \mu_1^r + 2\varepsilon_2) < \frac{(\lambda_1 - \mu_1)}{2}$, and $\frac{(\mu_1 - \lambda_1)}{2}s^r > 1 + 3\frac{\mu_1^r}{\mu_2^r}$, for all $r > r_3$. Then for $r > r_3$,

$$(5.23) \quad \mathcal{Q}^r(s^r) - \frac{\mu_1^r}{\mu_2^r}(\mathbf{C^r} - \mathbf{L^r} - 1) \leq 1 + 3\frac{\mu_1^r}{\mu_2^r} - \left(\frac{\mu_1 - \lambda_1}{2}\right)s^r < 0,$$

which, by definition of $\beta^r$, is a contradiction. This implies that for large enough $r$, on the set $\mathcal{Y}^r \cap \tilde{\mathcal{A}}^r$, $s^r \geq \beta^r$. Thus, on this set, for all $0 \leq s \leq \beta^r$ and large enough $r$,

$$\mathcal{Q}^r(s) = \mathcal{Q}^r(0) + \tilde{A}^r(s) - \tilde{S}^r(s)$$
$$\leq \mathcal{Q}^r(0) + \tilde{A}^r(s^r)$$
$$\leq \frac{\mu_1^r}{\mu_2^r}(\mathbf{C^r} - \mathbf{L^r} + 2) + 1 + (\lambda_1^r + \varepsilon_2)s^r$$
$$\leq \frac{3}{2}\frac{\mu_1^r}{\mu_2^r}(\mathbf{C^r} - \mathbf{L^r} + 1)$$
$$\leq \kappa''(\mathbf{C^r} - \mathbf{L^r} + 1),$$

where the last inequality follows on recalling that $\kappa'' > 2\frac{\mu_1}{\mu_2}$ and $2\frac{\mu_1}{\mu_2} > \frac{3\mu_1^r}{2\mu_2^r}$ for large $r$. This proves that $\mathbb{P}(\mathcal{Y}^r \cap \tilde{\mathcal{A}}^r) = 0$. Therefore, using (5.20), we have that each summand in the last term of (5.17) is bounded by $\mathbb{P}(\{\mathcal{Y}^r\}^c)$. Hence, using (5.22) and (5.17), we get, for $r$ sufficiently large (not depending on $t$),

$$\mathbb{P}(\mathcal{E}(r,t)) \leq (1 + 2n^r)\kappa_1 e^{-\kappa_2 r^2 t} + 2(n^r)^2\kappa_3 e^{-\kappa_4 s^r}$$
$$\leq \theta_1(r^4 t^2 + 1)(e^{-\theta_2 r^2 t} + r^{-\theta_3(c-1)\ell_0}),$$

for some positive constants $\theta_i$, $i = 1, 2, 3$. This completes the proof of Theorem 4.8. $\quad\square$

PROOF OF THEOREM 4.9. Let $d_0 \doteq d\ell_0$. We need to prove that, for sufficiently large $r$ (not depending on $t$),

$$(5.24) \qquad \mathbb{P}\left[\int_{[0,r^2 t)} \mathbf{I}_{\{Q_2^r(s) \geq d_0 \log r\}} dI_2^r(s) \neq 0\right]$$
$$\leq \gamma_1(1 + r^2 t)e^{-\gamma_2 r^2 t} + \gamma_3(1 + r^2 t)^2 r^{-\gamma_4 \ell_0}.$$

Fix $n \geq 1$. Note that, from (5.4), it follows that, for $s \in [\tau^r_{2n-1}, \tau^r_{2n})$, either $Q^r_3(s) \geq \mathbf{L^r} - 1 + \frac{\mu^r_1}{\mu^r_2} Q^r_1(s) > 0$, or $Q^r_3(s) \geq \mathbf{C^r} - 1 > 0$, for $r$ large enough. From the form of the control policy, we have that the idle-time process for the second server, $I^r_2(\cdot)$, does not increase during those intervals, and so the integrals over those intervals are zero. Thus, we need to consider only intervals of the form $[\tau^r_{2n-2}, \tau^r_{2n-1})$. We subdivide such intervals using a new sequence of stopping times as follows:

$$
\begin{aligned}
\tilde{\eta}^{r,n}_0 &\doteq \tau^r_{2n-2}, \\
\tilde{\eta}^{r,n}_{2m-1} &\doteq \min[\tau^r_{2n-1}, \inf\{s > \tilde{\eta}^{r,n}_{2m-2} | Q^r_3(s) \geq (\mathbf{C^r} - 1)\}], & m = 1, 2, \ldots, \\
\tilde{\eta}^{r,n}_{2m} &\doteq \min[\tau^r_{2n-1}, \inf\{s > \tilde{\eta}^{r,n}_{2m-1} | Q^r_3(s) < (\mathbf{C^r} - 1)\}], & m = 1, 2, \ldots, \\
\beta^{r,n}_m &\doteq \min\Big[\tilde{\eta}^{r,n}_{2m-1}, r^2 t, \\
& \qquad \inf\Big\{s > \tilde{\eta}^{r,n}_{2m-2} | Q^r_2(s) \geq \frac{d_0}{4} \log r\Big\}\Big], & m = 1, 2, \ldots.
\end{aligned}
\tag{5.25}
$$

Next, we estimate, exactly as in the proof of Theorem 4.8, how many such sub-intervals $[\tilde{\eta}^{r,n}_{2m-2}, \tilde{\eta}^{r,n}_{2m-1})$ there can be within $[0, r^2 t]$. Let $n^r$ be as in (5.7). Then from (5.8), we have that the probability in (5.24) is bounded by

$$
\begin{aligned}
\kappa_1 & e^{-\kappa_2 r^2 t} \\
& + \sum^{n^r}_{n=1} \mathbb{P}\Big[\int_{[\tau^r_{2n-2}, \tau^r_{2n-1} \wedge r^2 t)} \mathbf{I}_{\{Q^r_2(s) \geq d_0 \log r\}} \, dI^r_2(s) \neq 0, \tau^r_{2n-2} \leq r^2 t\Big].
\end{aligned}
\tag{5.26}
$$

Now within these intervals, consider the subintervals formed by $\tilde{\eta}^{r,n}_{\cdot}$'s. By the form of the policy, $Q^r_3(s) \geq \mathbf{C^r} - 1$, for $s \in [\tilde{\eta}^{r,n}_{2m-1}, \tilde{\eta}^{r,n}_{2m})$. Thus, $Q^r_3(s) = 0$ is possible only for $s \in [\tilde{\eta}^{r,n}_{2m-2}, \tilde{\eta}^{r,n}_{2m-1})$. Thus, we can conclude that

$$
\int_{[\tau^r_{2n-2}, \tau^r_{2n-1} \wedge r^2 t)} dI^r_2(s) = \sum^\infty_{m=1} \int_{[\tilde{\eta}^{r,n}_{2m-2}, \tilde{\eta}^{r,n}_{2m-1} \wedge r^2 t)} dI^r_2(s).
\tag{5.27}
$$

Next observe that on the set $\tilde{\eta}^{r,n}_{2m-2} < \tau^r_{2n-1} \wedge r^2 t$,

$$
Q^r_3(\tilde{\eta}^{r,n}_{2m-2}) \geq \mathbf{L^r} - 1 \qquad \text{for all } \tilde{\eta}^{r,n}_{2m-2} \geq \tilde{\eta}^{r,1}_2.
\tag{5.28}
$$

To see this, consider the case of $n \geq 2, m = 1$. Note that from (5.4), $Q^r_3(\tilde{\eta}^{r,n}_{2m-2}-) = Q^r_3(\tau^r_{2n-2}-) \geq \min\{\mathbf{L^r}, \mathbf{C^r} - 1\}$. This implies that $Q^r_3(\tilde{\eta}^{r,n}_0) \geq \mathbf{L^r} - 1$. For the case $n \geq 1, m \geq 2$, and $r$ sufficiently large, we have that $Q^r_3(\tilde{\eta}^{r,n}_{2m-2}-) \geq \mathbf{C^r} - 1 > \mathbf{L^r}$, which means $Q^r_3(\tilde{\eta}^{r,n}_{2m-2}) \geq \mathbf{L^r} - 1$. This proves (5.28).

Now define $m^r \doteq [r^2(\mu^r_2 + 1)] + 1$. Note that $\tilde{\eta}^{r,n}_{2m^r-1} < (\tau^r_{2n-1} \wedge r^2 t)$ implies that the queue-length $Q^r_3(\cdot)$ has crossed the threshold $\mathbf{C^r} - 1$ from below at least $m^r$ times before the time $r^2 t$, and each such up-crossing requires service of at least one job from $Buffer$ 2, implying $S^r_2(r^2 t) \geq m^r$. Using Corollary 5.3, we get

$$
\mathbb{P}(\tilde{\eta}^{r,n}_{2m^r-1} < \tau^r_{2n-1} \wedge r^2 t) \leq \mathbb{P}(S^r_2(r^2 t) \geq m^r) \leq \beta_3 e^{-\beta_4 r^2 t},
\tag{5.29}
$$

for all $r \geq \tilde{r}_5 \doteq r_1(\{\mu_2^r\}, 1)$, where $\beta_3 \doteq \varsigma_1(\{\mu_2^r\}, 1)$ , $\beta_4 \doteq \varsigma_2(\{\mu_2^r\}, 1)$, where $r_1(\cdot), \varsigma_1(\cdot), \varsigma_2(\cdot)$ are as in Corollary 5.3. Now, using (5.29) and (5.27), we write each summand of the second term in (5.26) as

$$\mathbb{P}\Big[\int_{[\tau_{2n-2}^r, \tau_{2n-1}^r \wedge r^2 t)} \mathbf{I}_{\{Q_2^r(s) \geq d_0 \log r\}} \, dI_2^r(s) \neq 0, \tau_{2n-2}^r \leq r^2 t\Big]$$

$$(5.30) \quad \leq \beta_3 e^{-\beta_4 r^2 t} + \sum_{m=1}^{m^r} \mathbb{P}\Big[\int_{[\tilde{\eta}_{2m-2}^{r,n}, \tilde{\eta}_{2m-1}^{r,n} \wedge r^2 t)} \mathbf{I}_{\{Q_2^r(s) \geq d_0 \log r\}} \, dI_2^r(s) \neq 0,$$
$$\tilde{\eta}_{2m-2}^{r,n} < \tau_{2n-1}^r \wedge r^2 t\Big].$$

From (5.26) and definitions of $n^r, m^r$, we get from (5.30) the following bound on the probability in (5.24):

$$\mathbb{P}\Big[\int_{[0, r^2 t)} \mathbf{I}_{\{Q_2^r(s) \geq d_0 \log r\}} \, dI_2^r(s) \neq 0\Big]$$

$$(5.31) \quad \leq \kappa_1 e^{-\kappa_2 r^2 t} + n^r \beta_3 e^{-\beta_4 r^2 t}$$
$$+ \sum_{n=1}^{n^r} \sum_{m=1}^{m^r} \mathbb{P}\Big(\int_{[\tilde{\eta}_{2m-2}^{r,n}, \tilde{\eta}_{2m-1}^{r,n} \wedge r^2 t)} \mathbf{I}_{\{Q_2^r(s) \geq d_0 \log r\}} \, dI_2^r(s) \neq 0,$$
$$\tilde{\eta}_{2m-2}^{r,n} < \tau_{2n-1}^r \wedge r^2 t\Big).$$

By definition of $\beta_m^{r,n}$, for $s \in [\tilde{\eta}_{2m-2}^{r,n}, \beta_m^{r,n})$, we have that $Q_2^r(s) < \frac{d_0}{4} \log r$. Therefore, for such $s$, the integrand on the right-hand side of (5.31) is zero. Hence, we have that

each summand in (5.31) is bounded by

$$\mathbb{P}\Big(\int_{[\beta_m^{r,n}, \tilde{\eta}_{2m-1}^{r,n} \wedge r^2 t)} \mathbf{I}_{\{Q_2^r(s) \geq d_0 \log r\}} \, dI_2^r(s) \neq 0, \beta_m^{r,n} < \tau_{2n-1}^r \wedge r^2 t\Big).$$

(5.32)

Now we make the following selections:

$$(5.33) \quad \begin{aligned} \varepsilon_1 &\in \Big(0, \min\Big\{\frac{(\mu_2 - \mu_3)}{8}, \frac{\mu_3}{8}, \frac{\lambda_2}{8}\Big\}\Big), \\ c &= \Big(1 + \frac{4}{\theta_3} + \frac{4(\mu_2 - \mu_3)}{\varsigma_2}\Big), \\ K &= 2 \max\{4, 16\lambda_2, 32\mu_2, 16\mu_3\}, \\ d &= c\Big(\frac{K}{(\mu_2 - \mu_3)/2}\Big), \\ c^* &= \frac{2d}{K}\ell_0, \\ \theta &= \frac{1}{2}\min\Big\{\frac{1}{4}, \frac{1}{32d}\Big\}, \end{aligned}$$

where $\theta_3$ is as in Theorem 4.8 and $\varsigma_2 \equiv \varsigma_2(\{\lambda_2^r\}, \varepsilon = \frac{\lambda_2}{2})$ is as in Corollary 5.3.

By Assumption 2.1 and choice of $\varepsilon_1$, we can find $r_1 \geq 1$, such that, for all $r \geq r_1$,

$$(5.34) \qquad \mu_2^r - \mu_3^r - 2\varepsilon_1 \geq \frac{\mu_2 - \mu_3}{2}, \qquad \mu_2^r + \varepsilon_1 < 2\mu_2, \qquad \mu_3^r + \varepsilon_1 < 2\mu_3.$$

Define, for $s > 0$, $\tilde{A}_2^r(s) \doteq A_2^r(\beta_m^{r,n} + s) - A_2^r(\beta_m^{r,n})$, $\tilde{S}_j^r(s) \doteq S_j^r(T_j^r(\beta_m^{r,n}) + s) - S_j^r(T_j^r(\beta_m^{r,n}))$, $j = 2, 3$. Now define

$$(5.35) \qquad \begin{aligned} \widetilde{\mathcal{A}_m^{r,n}} &\doteq \left\{ \sup_{\theta(c^* \log r) \leq s \leq (c^* \log r)} \left| \frac{\tilde{S}_j^r(s) - \mu_j^r s}{s} \right| \leq \varepsilon_1, \text{for } j = 2, 3 \right\}, \\ \mathcal{A}_m^{r,n} &\doteq \widetilde{\mathcal{A}_m^{r,n}} \cap \{\beta_m^{r,n} < \tau_{2n-1}^r \wedge r^2 t\}. \end{aligned}$$

Observing that $\{\beta_m^{r,n}\}$ are stopping times with respect to the filtration generated by the queue-length processes, using the strong Markov property of the Poisson processes, and using Corollary 5.2, we have, for some constant $C_1$ and function $C_2(\cdot)$ (not depending on $r$),

$$(5.36) \qquad \mathbb{P}[\{\widetilde{\mathcal{A}_m^{r,n}}\}^{\mathbf{c}}] \leq C_1 e^{-C_2(\varepsilon_1 \theta)(c^* \log r)}.$$

Using (5.36), we can write (5.32) as

$$\mathbb{P}\left( \int_{[\beta_m^{r,n}, \tilde{\eta}_{2m-1}^{r,n} \wedge r^2 t)} \mathbf{I}_{\{Q_2^r(s) \geq d_0 \log r\}} \, dI_2^r(s) \neq 0, \beta_m^{r,n} < \tau_{2n-1}^r \wedge r^2 t \right)$$

$$(5.37) \qquad \leq \mathbb{P}[\{\widetilde{\mathcal{A}_m^{r,n}}\}^{\mathbf{c}} \cap \{\beta_m^{r,n} < \tau_{2n-1}^r \wedge r^2 t\}]$$

$$\qquad\qquad + \mathbb{P}\left( \int_{[\beta_m^{r,n}, \tilde{\eta}_{2m-1}^{r,n} \wedge r^2 t)} \mathbf{I}_{\{Q_2^r(s) \geq d_0 \log r\}} \, dI_2^r(s) \neq 0, \mathcal{A}_m^{r,n} \right)$$

$$\qquad \leq C_1 e^{-C_2(\varepsilon_1 \theta)(c^* \log r)}$$

$$(5.38) \qquad\qquad + \mathbb{P}\bigg( \beta_m^{r,n} > \tilde{\eta}_{2m-2}^{r,n},$$

$$\qquad\qquad\qquad\qquad \int_{[\beta_m^{r,n}, \tilde{\eta}_{2m-1}^{r,n} \wedge r^2 t)} \mathbf{I}_{\{Q_2^r(s) \geq d_0 \log r\}} \, dI_2^r(s) \neq 0, \mathcal{A}_m^{r,n} \bigg)$$

$$(5.39) \qquad\qquad + \mathbb{P}\bigg( \beta_m^{r,n} = \tilde{\eta}_{2m-2}^{r,n},$$

$$\qquad\qquad\qquad\qquad \int_{[\beta_m^{r,n}, \tilde{\eta}_{2m-1}^{r,n} \wedge r^2 t)} \mathbf{I}_{\{Q_2^r(s) \geq d_0 \log r\}} \, dI_2^r(s) \neq 0, \mathcal{A}_m^{r,n} \bigg).$$

Now consider the event corresponding to the probability in (5.38):

$$(5.40) \qquad \left\{ \beta_m^{r,n} > \tilde{\eta}_{2m-2}^{r,n}, \int_{[\beta_m^{r,n}, \tilde{\eta}_{2m-1}^{r,n} \wedge r^2 t)} \mathbf{I}_{\{Q_2^r(s) \geq d_0 \log r\}} \, dI_2^r(s) \neq 0, \mathcal{A}_m^{r,n} \right\} \doteq \mathcal{B}_m^{r,n}.$$

We claim that, for large values of $r$,

$$(5.41) \quad \mathbb{P}\Big[ \sup_{0 \le s \le (\theta c^* \log r)} \Big| Q_2^r(\beta_m^{r,n} + s) - \Big(\frac{d_0}{4} \log r\Big) \Big| \ge \frac{1}{2}\Big(\frac{d_0}{4} \log r\Big); \mathcal{B}_m^{r,n} \Big]$$
$$\le \varsigma_1' r^{-\varsigma_2' \theta c^*},$$

where $\varsigma_1' \equiv \varsigma_1(\{\mu_2^r\}, \varepsilon = \frac{\mu_2}{2} \wedge 1) + \varsigma_1(\{\lambda_2^r\}, \varepsilon = \frac{\lambda_2}{2} \wedge 1)$ and $\varsigma_2' \equiv \min\{\varsigma_2(\{\mu_2^r\}, \varepsilon = \frac{\mu_2}{2} \wedge 1), \varsigma_2(\{\lambda_2^r\}, \varepsilon = \frac{\lambda_2}{2} \wedge 1)\}$, where $\varsigma_1$ and $\varsigma_2$ are as in Corollary 5.3. To see this, note that, on $\mathcal{B}_m^{r,n}$, by definition of $\beta_m^{r,n}$ (recall that on the set $\mathcal{B}_m^{r,n}$, $\beta_m^{r,n} < r^2 t$), we have $Q_2^r(\beta_m^{r,n}) \ge \frac{d_0}{4} \log r$. And in order for $Q_2^r(\beta_m^{r,n} + s)$ to decrease by $\frac{1}{2}(\frac{d_0}{4} \log r)$ inside $[0, (\theta c^* \log r)]$, we need the number of services from *Buffer* 2 (by *Server* 1) to be greater than or equal to $\frac{1}{2}(\frac{d_0}{4} \log r)$ in $(\theta c^* \log r)$ time. On the other hand, since on the set $\mathcal{B}_m^{r,n}$, using $\beta_m^{r,n} > \tilde{\eta}_{2m-2}^{r,n}$, we have $Q_2^r(\beta_m^{r,n-}) < \frac{d_0}{4} \log r$, we can conclude that $Q_2^r(\beta_m^{r,n})$ can be at most $\frac{d_0}{4} \log r + 1$. And $Q_2^r(\beta_m^{r,n} + \cdot) \ge \frac{3}{2}(\frac{d_0}{4} \log r)$ inside $[0, (\theta c^* \log r)]$ implies the number of arrivals in *Buffer* 2 needs to be greater than or equal to $\frac{1}{2}(\frac{d_0}{4} \log r) - 1 > \frac{1}{4}(\frac{d_0}{4} \log r)$ in time $(\theta c^* \log r)$, for $r$ large. Thus, for $r$ sufficiently large, we can bound the probability in (5.41) by

$$(5.42) \quad \mathbb{P}\Big[ \tilde{S}_2^r(\theta c^* \log r) \ge \frac{1}{2}\Big(\frac{d_0}{4}\Big) \log r \Big] + \mathbb{P}\Big[ \tilde{A}_2^r(\theta c^* \log r) \ge \frac{1}{4}\Big(\frac{d_0}{4}\Big) \log r \Big].$$

Also note that, if $\varepsilon_2 = \min\{1, \frac{\mu_2}{2}\}$ and $\varepsilon_2' = \min\{1, \frac{\lambda_2}{2}\}$, by the choices made in (5.33), we have, for $r$ large enough,

$$\theta c^*(\mu_2^r + \varepsilon_2) < \theta \frac{2d_0}{K}(2\mu_2) < \frac{2\mu_2 d_0}{K} < \frac{1}{2}\Big(\frac{d_0}{4}\Big),$$
$$\theta c^*(\lambda_2^r + \varepsilon_2') < \theta \frac{2d_0}{K}(2\lambda_2) < \frac{2\lambda_2 d_0}{K} < \frac{1}{4}\Big(\frac{d_0}{4}\Big).$$

Using this observation and Corollary 5.3, we get that the sum in (5.42) is bounded by

$$\mathbb{P}[\tilde{S}_2^r(\theta c^* \log r) > (\mu_2^r + \varepsilon_2)\theta c^* \log r] + \mathbb{P}[\tilde{A}_2^r(\theta c^* \log r) > (\lambda_2^r + \varepsilon_2')\theta c^* \log r]$$
$$\le \varsigma_1' r^{-\varsigma_2' \theta c^*}.$$

This completes the proof of the claim in (5.41).

Now, using (5.41), for large values of $r$, we can bound (5.38) by

$$\mathbb{P}\Big( \beta_m^{r,n} > \tilde{\eta}_{2m-2}^{r,n}, \int_{[\beta_m^{r,n}, \tilde{\eta}_{2m-1}^{r,n} \wedge r^2 t]} \mathbf{I}_{\{Q_2^r(s) \ge d_0 \log r\}} \, dI_2^r(s) \ne 0, \mathcal{A}_m^{r,n} \Big)$$
$$\le \varsigma_1' r^{-\varsigma_2' \theta c^*}$$
$$+ \mathbb{P}\Big[ \tilde{\eta}_{2m-1}^{r,n} - \beta_m^{r,n} - \theta c^* \log r \le \frac{d_0}{K} \log r,$$

(5.43)
$$\sup_{0 \le s \le \theta c^* \log r} Q_2^r(\beta_m^{r,n} + s) \le \frac{3}{2}\Big(\frac{d_0}{4}\log r\Big), \mathcal{B}_m^{r,n}\Big]$$

$$+ \mathbb{P}\Big[\tilde{\eta}_{2m-1}^{r,n} - \beta_m^{r,n} - \theta c^* \log r > \frac{d_0}{K}\log r,$$

(5.44)
$$\inf_{0 \le s \le \theta c^* \log r} Q_2^r(\beta_m^{r,n} + s) \ge \frac{1}{2}\Big(\frac{d_0}{4}\log r\Big), \mathcal{B}_m^{r,n}\Big].$$

Now we get a bound on each of (5.43) and (5.44). For the event in (5.43), note that $Q_2^r(\beta_m^{r,n} + s) \le \frac{3}{2}\big(\frac{d_0}{4}\log r\big)$ for $s \le \theta c^* \log r$, and within an additional $\frac{d_0}{K}\log r$ units of time, $Q_2^r(s)$ becomes greater than or equal to $d_0 \log r$ [see the definition of the set $\mathcal{B}_m^{r,n}$ in (5.40)]. This implies that there are more than $\big(\frac{d_0}{4}\log r\big)$ arrivals in *Buffer* 2 in time $\frac{d_0}{K}\log r$. Recalling the definition of $\varepsilon_2' = \min\{1, \lambda_2/2\}$ and the choices made in (5.33), we have that, for $r$ large enough,

$$(\lambda_2^r + \varepsilon_2')\frac{d_0}{K} < (2\lambda_2)\frac{d_0}{K} < \frac{d_0}{4}.$$

Observing that $\{\beta_m^{r,n}\}$ are stopping times with respect to the filtration generated by the queue-length processes, using the strong Markov property of the Poisson processes, we get that the distribution (conditioned on the $\sigma$-field generated by queue-length processes stopped at $\beta_m^{r,n}$) of $\tilde{A}_2^r(\theta c^* \log r + \cdot) - \tilde{A}_2^r(\theta c^* \log r)$ is the same as that of $A_2^r(\cdot)$. This, together with the display above, yields that, for $r$ sufficiently large, (5.43) is bounded by

$$\mathbb{P}\Big[A_2^r\Big(\frac{d_0}{K}\log r\Big) > \Big(\frac{d_0}{4}\Big)\log r\Big] \le \mathbb{P}\Big[A_2^r\Big(\frac{d_0}{K}\log r\Big) > (\lambda_2^r + \varepsilon_2')\frac{d_0}{K}\log r\Big]$$
$$\le \varsigma_1'' r^{-\varsigma_2'' d_0/K},$$

where $\varsigma_i'' \equiv \varsigma_i(\{\lambda_2^r\}, \varepsilon = \frac{\lambda_2}{2} \wedge 1)\}$, $i = 1, 2$, are as in Corollary 5.3.

Next we show that, for $r$ sufficiently large, (5.44) is zero. Note that by the choices made in (5.33), on $\mathcal{A}_m^{r,n}$, we have that, for $r$ sufficiently large,

(5.45)
$$\tilde{S}_2^r(s) \le (\mu_2^r + \varepsilon_1)s \le (2\mu_2)\frac{2d_0}{K}\log r < \frac{1}{2}\Big(\frac{d_0}{4}\log r\Big)$$
$$\text{for all } s \in \Big[\theta c^* \log r, \frac{2d_0}{K}\log r\Big].$$

Since on the set in (5.44) $Q_2^r(\beta_m^{r,n} + \theta c^* \log r) \ge \frac{1}{2}\big(\frac{d_0}{4}\log r\big)$, this means $Q_2^r(\beta_m^{r,n} + s)$ never becomes zero for $s$ in the interval $[\theta c^* \log r, (2d_0/K)\log r]$. So, on the set in (5.45), $Q_2^r(\beta_m^{r,n} + s)$ never becomes zero for $s$ in the interval $[0, (2d_0/K)\log r]$. Hence, using the fact that our policy requires *Server* 1 to work on *Buffer* 2 continuously in the interval $[\beta_m^{r,n}, \beta_m^{r,n} + \theta c^* \log r + \frac{d_0}{K} \log r]$, on the set in (5.44), we have that

$$
\begin{aligned}
(5.46) \quad Q_3^r &\left(\beta_m^{r,n} + \theta c^* \log r + \frac{d_0}{K} \log r\right) \\
&\geq \tilde{S}_2^r\left(\theta c^* \log r + \frac{d_0}{K} \log r\right) - \tilde{S}_3^r\left(\theta c^* \log r + \frac{d_0}{K} \log r\right) \\
&\geq (\mu_2^r - \mu_3^r - 2\varepsilon_1)\left(\theta c^* \log r + \frac{d_0}{K} \log r\right) \\
&\geq \frac{(\mu_2 - \mu_3)}{2}\frac{d_0}{K}\log r \\
&= c_0 \log r \geq \mathbf{C^r}.
\end{aligned}
$$

However, (5.46) is a contradiction to the fact that, on the set in (5.44), we have $\tilde{\eta}_{2m-1}^{r,n} > \beta_m^{r,n} + \theta c^* \log r + \frac{d_0}{K} \log r$. This proves that (5.44) is zero. Thus, the term (5.38) is bounded by

$$
(5.47) \qquad\qquad \varsigma_1' r^{-\varsigma_2' \theta c^*} + \varsigma_1'' r^{-\varsigma_2'' d_0/K}.
$$

Now we consider (5.39). First note that for $n = 1, m = 1$, (5.39) is zero, since $Q_2^r(0) = 0$. For all other $n \geq 1, m \geq 1$, consider the event corresponding to the probability in (5.39):

$$
\left\{\beta_m^{r,n} = \tilde{\eta}_{2m-2}^{r,n}, \int_{[\beta_m^{r,n}, \tilde{\eta}_{2m-1}^{r,n} \wedge r^2 t)} \mathbf{I}_{\{Q_2^r(s) \geq d_0 \log r\}} \, dI_2^r(s) \neq 0, \mathcal{A}_m^{r,n}\right\} \doteq \mathcal{C}_m^{r,n}.
$$

We claim that, for large values of $r$,

$$
\begin{aligned}
(5.48) \quad \mathbb{P}\bigg(&\bigg[\Big\{\inf_{0 \leq s \leq \theta c^* \log r} Q_3^r(\beta_m^{r,n} + s) \leq \frac{1}{2}\ell_0 \log r\Big\} \\
&\cup \Big\{\inf_{0 \leq s \leq (\theta c^* \log r)} Q_2^r(\beta_m^{r,n} + s) \leq \frac{1}{2}\Big(\frac{d_0}{4}\log r\Big)\Big\}\bigg], \mathcal{C}_m^{r,n}\bigg) \\
&\leq \rho_1 r^{-\rho_2 \theta c^*},
\end{aligned}
$$

where $\rho_1 \doteq \varsigma_1(\{\mu_3^r\}, \varepsilon = \frac{\mu_3}{2} \wedge 1) + \varsigma_1(\{\mu_2^r\}, \varepsilon = \frac{\mu_2}{2} \wedge 1)$ and $\rho_2 \doteq \min\{\varsigma_2(\{\mu_3^r\}, \varepsilon = \frac{\mu_3}{2} \wedge 1), \varsigma_2(\{\mu_2^r\}, \varepsilon = \frac{\mu_2}{2} \wedge 1)\}$ and $\varsigma_i, i = 1, 2$, are as in Corollary 5.3.

To see the claim, note that, on $\mathcal{C}_m^{r,n}$, $\beta_m^{r,n} = \tilde{\eta}_{2m-2}^{r,n}$ and from (5.28), $Q_3^r(\tilde{\eta}_{2m-2}^{r,n}) \geq \ell_0 \log r - 1$. And if $r$ is large enough so that $\frac{1}{2}(\ell_0 \log r) - 1 > \frac{\ell_0}{4} \log r$, then in order for $Q_3^r(\beta_m^{r,n} + \cdot)$ to decrease by $\frac{1}{2}(\ell_0 \log r) - 1$ in $(\theta c^* \log r)$ time, we need the number of Class 3 services in that time interval to be greater than or equal to $(\frac{\ell_0}{4} \log r)$. On the other hand, note that, by definition of $\beta_m^{r,n}$, we have $Q_2^r(\beta_m^{r,n}) \geq \frac{d_0}{4} \log r$. And in order for $Q_2^r(\beta_m^{r,n} + \cdot)$ to decrease by $\frac{1}{2}(\frac{d_0}{4} \log r)$ inside $[0, (\theta c^* \log r)]$, we need the number of Class 2 services in that time interval to be greater than or equal to $\frac{1}{2}(\frac{d_0}{4} \log r)$. So for large values of $r$, the probability in (5.48) is bounded by

$$
(5.49) \quad \mathbb{P}\bigg[\tilde{S}_3^r(\theta c^* \log r) \geq \frac{\ell_0}{4} \log r\bigg] + \mathbb{P}\bigg[\tilde{S}_2^r(\theta c^* \log r) \geq \frac{1}{2}\Big(\frac{d_0}{4}\Big) \log r\bigg].
$$

Also note that, by the choices made in (5.33), if we set $\varepsilon_3 \doteq \min\{1, \mu_3/2\}, \varepsilon_4 \doteq \min\{1, \mu_2/2\}$, we have, for $r$ large enough, that

$$\theta c^*(\mu_3^r + \varepsilon_3) < \theta \frac{2d_0}{K}(2\mu_3) = \theta d \frac{4\mu_3}{K}\ell_0 < \frac{\ell_0}{4},$$

$$\theta c^*(\mu_2^r + \varepsilon_4) < \theta \frac{2d_0}{K}(2\mu_2) < \frac{2\mu_2 d_0}{K} < \frac{1}{2}\left(\frac{d_0}{4}\right).$$

Using the above observations and Corollary 5.3, we get that, for large $r$, the sum in (5.49) is bounded by

$$\mathbb{P}[\tilde{S}_3^r(\theta c^* \log r) > (\mu_3^r + \varepsilon_3)\theta c^* \log r]$$
$$+ \mathbb{P}[\tilde{S}_2^r(\theta c^* \log r) > (\mu_2^r + \varepsilon_4)\theta c^* \log r] \leq \rho_1 r^{-\rho_2 \theta c^*}.$$

This completes the proof of the claim (5.48). Now using (5.48), for large values of $r$, we can bound (5.39) as follows:

$$\mathbb{P}\left(\beta_m^{r,n} = \tilde{\eta}_{2m-2}^{r,n}, \int_{[\beta_m^{r,n}, \tilde{\eta}_{2m-1}^{r,n} \wedge r^2 t)} \mathbf{I}_{\{Q_2^r(s) \geq d_0 \log r\}} \, dI_2^r(s) \neq 0, \mathcal{A}_m^{r,n}\right)$$

$$\leq \rho_1 r^{-\rho_2 \theta c^*}$$

$$+ \mathbb{P}\left[\tilde{\eta}_{2m-1}^{r,n} - \beta_m^{r,n} - \theta c^* \log r \leq \frac{2c_0}{(\mu_2 - \mu_3)} \log r,\right.$$

(5.50) $$Q_2^r(\beta_m^{r,n} + s) > \frac{1}{2}\left(\frac{d_0}{4} \log r\right) \text{ for } s \leq \theta c^* \log r,$$

$$\left.\inf_{0 \leq s \leq \theta c^* \log r} Q_3^r(\beta_m^{r,n} + s) \geq \frac{\ell_0}{2} \log r, \mathcal{C}_m^{r,n}\right]$$

$$+ \mathbb{P}\left[\tilde{\eta}_{2m-1}^{r,n} - \beta_m^{r,n} - \theta c^* \log r > \frac{2c_0}{(\mu_2 - \mu_3)} \log r,\right.$$

(5.51) $$Q_2^r(\beta_m^{r,n} + s) > \frac{1}{2}\left(\frac{d_0}{4} \log r\right) \text{ for } s \leq \theta c^* \log r,$$

$$\left.\inf_{0 \leq s \leq \theta c^* \log r} Q_3^r(\beta_m^{r,n} + s) \geq \frac{\ell_0}{2} \log r, \mathcal{C}_m^{r,n}\right].$$

We will next show that both the terms (5.50) and (5.51) are zero. First observe that, from (5.45), we have that, on $\mathcal{A}_m^{r,n}$, for $r$ sufficiently large, $Q_2^r(\beta_m^{r,n} + \cdot)$ never becomes zero in the interval $[\theta c^* \log r, (2d_0/K) \log r]$. Thus, on the sets corresponding to terms (5.50) and (5.51), $Q_2^r(\beta_m^{r,n} + \cdot)$ is not zero on $[0, (2d_0/K) \log r]$.

On the event in (5.50), we have

$$Q_3^r(\beta_m^{r,n} + \theta c^* \log r) \geq \frac{\ell_0}{2} \log r$$

and

$$[(\beta_m^{r,n} + \theta c^* \log r), \tilde{\eta}_{2m-1}^{r,n}]$$
$$\subseteq \left[ (\beta_m^{r,n} + \theta c^* \log r), (\beta_m^{r,n} + \theta c^* \log r) + \frac{2c_0}{(\mu_2 - \mu_3)} \right] \log r \Big].$$

Thus, by definition of $\mathcal{C}_m^{r,n}$ and conditions of the event in (5.50), we must have that $Q_3^r(\beta_m^{r,n} + \theta c^* \log r + s)$ is zero for some $s$ in $[0, \frac{2c_0}{(\mu_2 - \mu_3)} \log r]$. This means that, for some $s$ in the above interval, $[Q_3^r(\beta_m^{r,n} + \theta c^* \log r) - Q_3^r(\beta_m^{r,n} + \theta c^* \log r + s)] \geq \frac{\ell_0}{4} \log r$. Now since $Q_2^r(\beta_m^{r,n} + \cdot)$ never becomes zero in the interval $[0, (2d_0/K) \log r]$, this decrease in $Q_3^r(s)$ is bounded by $[\tilde{S}_3^r((\theta c^* \log r) + s) - \tilde{S}_2^r((\theta c^* \log r) + s)] - [\tilde{S}_3^r(\theta c^* \log r) - \tilde{S}_2^r(\theta c^* \log r)]$. Hence, the probability in (5.50) is bounded above by

$$\mathbb{P}\Bigg[ \sup_{0 \leq s \leq 2c_0 \log r/(\mu_2 - \mu_3)} (\tilde{S}_3^r(\theta c^* \log r + s)$$
$$- \tilde{S}_2^r(\theta c^* \log r + s)) - (\tilde{S}_3^r(\theta c^* \log r) - \tilde{S}_2^r(\theta c^* \log r)) \geq \frac{\ell_0}{4} \log r; \mathcal{A}_m^{r,n} \Bigg].$$

We claim that, for all $r$ large enough, the above probability is zero. To see the claim, note that, for all $s \leq \frac{2c_0 \log r}{(\mu_2 - \mu_3)}$, we have from (5.33) that $\theta c^* \log r + s \in [\theta c^* \log r, c^* \log r]$. Thus, by definition of $\mathcal{A}_m^{r,n}$ and (5.33), we get, for all such $s$, that

$$(\tilde{S}_3^r(\theta c^* \log r + s) - \tilde{S}_2^r(\theta c^* \log r + s)) - (\tilde{S}_3^r(\theta c^* \log r) - \tilde{S}_2^r(\theta c^* \log r))$$
$$< (\mu_3^r - \mu_2^r + 2\varepsilon_1)(\theta c^* \log r + s) - (\mu_3^r - \mu_2^r - 2\varepsilon_1)(\theta c^* \log r)$$
$$= 4\varepsilon_1(\theta c^* \log r) + (\mu_3^r - \mu_2^r + 2\varepsilon_1)s$$
$$\leq 4\varepsilon_1(\theta c^* \log r) - \tfrac{1}{2}(\mu_2 - \mu_3)s$$
$$\leq (4\varepsilon_1 - \tfrac{1}{2}(\mu_2 - \mu_3))(\theta c^* \log r)$$
$$< 0,$$

for $r$ large enough so that $(\mu_3^r - \mu_2^r + 2\varepsilon_1) < -(\mu_2 - \mu_3)/2$. This proves the claim. Thus, the expression in (5.50) is zero. We now show that (5.51) is zero as well. To see that, recall that, on the event in (5.51), $Q_2^r(\beta_m^{r,n} + \cdot)$ never becomes zero in the interval $[0, (2d_0/K) \log r]$. This implies that

$$(5.52) \quad \begin{aligned} Q_3^r &\left( \beta_m^{r,n} + \theta c^* \log r + \frac{2c_0 \log r}{(\mu_2 - \mu_3)} \right) \\ &\geq \tilde{S}_2^r \left( \theta c^* \log r + \frac{2c_0 \log r}{(\mu_2 - \mu_3)} \right) - \tilde{S}_3^r \left( \theta c^* \log r + \frac{2c_0 \log r}{(\mu_2 - \mu_3)} \right) \\ &\geq (\mu_2^r - \mu_3^r - 2\varepsilon_1) \left( \theta c^* \log r + \frac{2c_0 \log r}{(\mu_2 - \mu_3)} \right) \\ &\geq \frac{(\mu_2 - \mu_3)}{2} \frac{2c_0 \log r}{(\mu_2 - \mu_3)} = c_0 \log r \geq \mathbf{C^r}. \end{aligned}$$

However, (5.52) contradicts the definition of $\tilde{\eta}^{r,n}_{2m-1}$ in view of the fact that, on the set in (5.51), $\tilde{\eta}^{r,n}_{2m-1} > \beta^{r,n}_m + \theta c^* \log r + \frac{2c_0 \log r}{(\mu_2 - \mu_3)}$. This proves that (5.51) is zero. Hence, we have proved that (5.39) is bounded by $\rho_1 r^{-\rho_2 \theta c^*}$. Using the above observation and (5.47), the term in (5.37) is bounded by

$$\Lambda(r) \doteq C_1 e^{-[C_2(\varepsilon_1 \theta)](c^* \log r)} + \varsigma'_1 r^{-\varsigma'_2 \theta c^*} + \varsigma''_1 r^{-\varsigma''_2 d_0/K} + \rho_1 r^{-\rho_2 \theta c^*}.$$

Using (5.32) and (5.31), we get the following bound on the left-hand side of (5.24), for large enough $r$:

$$
\begin{aligned}
(5.53) \qquad &\mathbb{P}\Big[\int_{[0,t)} \mathbf{I}_{\{\hat{Q}^r_2(s) \geq d_0 \log r/r\}}\, d\hat{I}^r_2(s) \neq 0\Big] \\
&\leq \kappa_1 e^{-\kappa_2 r^2 t} + n^r \beta_3 e^{-\beta_4 r^2 t} + n^r m^r \Lambda(r) \\
&\leq \gamma_1 (1 + r^2 t) e^{-\gamma_2 r^2 t} + \gamma_3 (1 + r^2 t)^2 r^{-\gamma_4 \ell_0},
\end{aligned}
$$

for some constants $\gamma_i > 0, i = 1, \dots, 4$, which are independent of $t$. This completes the proof of the theorem. $\square$

Proof of Theorem 4.11. We will only prove (4.41). The proof of (4.40) is similar and therefore is omitted. Consider first the case $i = 1$. In view of (4.45), the main step is to obtain bounds for the following integrals:

$$(5.54) \qquad \int_T^\infty e^{-\gamma t} \mathbb{E}\Big[\Big\{\sup_{0 \leq s \leq t} |\hat{X}^r_i(s)|\Big\}^2\Big]\, dt, \qquad i = 1, 2.$$

By definition of $\hat{X}^r_i(\cdot)$, we have that, for $i = 1, 2$,

$$
\begin{aligned}
(5.55) \qquad \sup_{0 \leq s \leq t} (\hat{X}^r_i(s))^2 &\leq 3\Big(\sup_{0 \leq s \leq t} |\hat{A}^r_i(s)|\Big)^2 \\
&\quad + 3\Big(\sup_{0 \leq s \leq t} |\hat{S}^r_i(s)|\Big)^2 + 3\Big(r\mu^r_i\Big(\frac{\lambda^r_i}{\mu^r_i} - \frac{\lambda_i}{\mu_i}\Big)t\Big)^2.
\end{aligned}
$$

By Doob's maximal inequality [for the martingale $(A^r_i(s) - \lambda^r_i s)$], we have

$$
\begin{aligned}
(5.56) \qquad &\int_T^\infty e^{-\gamma t} \mathbb{E}\Big[\Big\{\sup_{0 \leq s \leq t} |\hat{A}^r_i(s)|\Big\}^2\Big]\, dt \\
&= r^{-2} \int_T^\infty e^{-\gamma t} \mathbb{E}\Big[\Big\{\sup_{0 \leq s \leq r^2 t} |A^r_i(s) - \lambda^r_i s|\Big\}^2\Big]\, dt \\
&\leq 4\lambda^r_i \int_T^\infty t e^{-\gamma t}\, dt.
\end{aligned}
$$

In a similar way, one shows that

$$(5.57) \qquad \int_T^\infty e^{-\gamma t} \mathbb{E}\Big[\Big\{\sup_{0 \leq s \leq t} |\hat{S}^r_i(s)|\Big\}^2\Big]\, dt \leq 4\mu^r_i \int_T^\infty t e^{-\gamma t}\, dt.$$

Combining (5.56), (5.57), (5.55) and using Assumption 2.2, we obtain

$$(5.58) \quad \limsup_{T \to \infty} \limsup_{r \to \infty} \int_T^\infty e^{-\gamma t} \mathbb{E}\left[ \left\{ \sup_{0 \le s \le t} |\hat{X}_i^r(s)| \right\}^2 \right] dt = 0, \qquad i = 1, 2.$$

Finally, combining (5.58) and (4.45) with the fact that $\Gamma(\cdot)$ is Lipschitz continuous, we have (4.41) for $i = 1$.

Proof of (4.41), for $i = 2$, is similar. We will only prove the key steps. Let

$$(5.59) \qquad\qquad Y_t^r \doteq \int_{[0,t)} \mathbf{I}_{\{\hat{Q}_2^r(s) \ge d_0 \log r/r\}} \, d\hat{I}_2^r(s).$$

From (4.49), it is clear that we need to get an estimate on

$$(5.60) \quad \int_T^\infty e^{-\gamma t} \mathbb{E}(\{Y_t^r\}^2) \, dt = \int_T^\infty e^{-\gamma t} \int_0^\infty \mathbb{P}(Y_t^r > \sqrt{u}) \, du \, dt.$$

Theorem 4.9 and the fact that $\hat{I}_2^r(s) \le rs$ yields the following bound on the integral in (5.60):

$$
\begin{aligned}
\int_0^\infty \mathbb{P}(Y_t^r > \sqrt{u}) \, du &= \int_0^{r^2 t^2} \mathbb{P}\left( \int_{[0,t)} \mathbf{I}_{\{\hat{Q}_2^r(s) \ge d_0 \log r/r\}} \, d\hat{I}_2^r(s) > \sqrt{u} \right) du \\
(5.61) \qquad &\le r^2 t^2 \mathbb{P}\left( \int_{[0,t)} \mathbf{I}_{\{\hat{Q}_2^r(s) \ge d_0 \log r/r\}} \, d\hat{I}_2^r(s) \ne 0 \right) \\
&\le r^2 t^2 (\gamma_1 (1 + r^2 t) e^{-\gamma_2 r^2 t} + \gamma_3 (1 + r^2 t)^2 r^{-\gamma_4 \ell_0}).
\end{aligned}
$$

Substituting the above estimate in (5.60), one obtains after some straightforward calculations that

$$(5.62) \qquad\qquad \limsup_{T \to \infty} \limsup_{r \to \infty} \int_T^\infty e^{-\gamma t} \mathbb{E}(Y_t^r)^2 \, dt = 0.$$

Using (4.54), it follows that

$$(5.63) \quad \limsup_{T \to \infty} \limsup_{r \to \infty} \int_T^\infty e^{-\gamma t} \mathbb{E}\left( \sup_{0 \le s \le t} \frac{\hat{Q}_2^r(s)}{\mu_3^r} \mathbf{I}_{\{\hat{Q}_2^r(s) < 2 d_0 \log r/r\}} \right)^2 dt = 0.$$

Now, as in the first half of the proof [see (5.58)], we can prove that

$$(5.64) \quad \limsup_{T \to \infty} \limsup_{r \to \infty} \int_T^\infty e^{-\gamma t} \mathbb{E}\left[ \left\{ \sup_{0 \le s \le t} |\hat{X}_i^r(s)| \right\}^2 \right] dt = 0, \qquad i = 2, 3.$$

Thus, from (4.49), (5.62), (5.63) and the Lipschitz property of the Skorohod map, we get (4.41) for $i = 2$. This completes the proof of the theorem. $\quad\square$

## APPENDIX: PROOFS OF LEMMAS 4.5 AND 4.6

PROOF OF LEMMA 4.5. From (3.3), it follows that

$$(A.1) \qquad (\bar{A}^r(\cdot), \bar{S}^r(\cdot)) \Rightarrow (\lambda(\cdot), \mu(\cdot)) \qquad \text{as } r \to \infty,$$

where $\lambda(t) = \lambda t, \mu(\cdot) = \mu t; t \geq 0$. Also, it follows from (2.7) and definition of the fluid-scaled processes in (2.8) that $\bar{T}^r(\cdot)$ is uniformly Lipschitz continuous with Lipschitz constant less than or equal to 1. This fact and (A.1) imply that

$$(A.2) \qquad \{\bar{A}^r(\cdot), \bar{S}^r(\cdot), \bar{T}^r(\cdot)\}_{r \in \mathbb{S}} \text{ is } \mathcal{C}\text{-tight.}$$

Now, by definition of the queue-length process (2.1) and fluid-scaled processes (2.8), we have

$$(A.3) \qquad \begin{aligned} \bar{Q}_1^r(t) &= \bar{A}_1^r(t) - \bar{S}_1^r(\bar{T}_1^r(t)), \\ \bar{Q}_2^r(t) &= \bar{A}_2^r(t) - \bar{S}_2^r(\bar{T}_2^r(t)), \\ \bar{Q}_3^r(t) &= \bar{S}_2^r(\bar{T}_2^r(t)) - \bar{S}_3^r(\bar{T}_3^r(t)), \\ \bar{I}_1^r(t) &= t - \bar{T}_1^r(t) - \bar{T}_2^r(t), \\ \bar{I}_2^r(t) &= t - \bar{T}_3^r(t). \end{aligned}$$

From (A.2), (A.3) and Lemma 3.14.1 of [2], we get

$$(A.4) \qquad \{\bar{Q}^r(\cdot), \bar{I}^r(\cdot)\}_{r \in \mathbb{S}} \text{ is } \mathcal{C}\text{-tight.}$$

Combining (A.2) and (A.4), we have (4.4). $\square$

PROOF OF LEMMA 4.6. From Lemma 4.5, we have

$$(A.5) \qquad \{\bar{Q}^{r'}(\cdot), \bar{A}^{r'}(\cdot), \bar{S}^{r'}(\cdot), \bar{T}^{r'}(\cdot), \bar{I}^{r'}(\cdot)\}_{r' \geq 1} \text{ is } \mathcal{C}\text{-tight.}$$

Thus, it is sufficient to show that all weak limit-points of the above sequence are given by the right-hand side of (4.7).

Suppose that $(\bar{Q}(\cdot), \bar{A}(\cdot), \bar{S}(\cdot), \bar{T}(\cdot), \bar{I}(\cdot))$ is a limit-point of the sequence in (A.5), obtained along a subsequence indexed by $r''$. Using the Skorohod representation theorem, we can assume that this convergence takes place almost surely, uniformly on compacts:

$$(A.6) \qquad \begin{aligned} (\bar{Q}^{r''}(\cdot), \bar{A}^{r''}(\cdot), \bar{S}^{r''}(\cdot), \bar{T}^{r''}(\cdot), \bar{I}^{r''}(\cdot)) \\ \to (\bar{Q}(\cdot), \bar{A}(\cdot), \bar{S}(\cdot), \bar{T}(\cdot), \bar{I}(\cdot)) \qquad \text{as } r'' \to \infty. \end{aligned}$$

From (A.1), we have that $\bar{A}(\cdot) = \lambda(\cdot)$ and $\bar{S}(\cdot) = \mu(\cdot)$. Recall that, by assumption, $\lim_{r'' \to \infty} \hat{J}^{r''}(T^{r''}) = \underline{J}(\{T^{r''}\}) < \infty$. Thus, using Fatou's lemma, we get

$$0 = \lim_{r'' \to \infty} \frac{1}{r''} \hat{J}^{r''}(T^{r''}) \geq \mathbf{E}\left( \int_0^\infty e^{-\gamma t} \liminf_{r'' \to \infty} (h \cdot \bar{Q}^{r''}(t)) \, dt \right)$$

$$= \mathbf{E}\left( \int_0^\infty e^{-\gamma t} (h \cdot \bar{Q}(t)) \, dt \right).$$

Since $h_i > 0, i = 1, 2, 3$, and $\bar{Q}$ has continuous paths, a.s., we have from the above equation that $\bar{Q}(\cdot) \equiv 0$. Using this along with (A.3) and (A.1), we now see that, for all $t \geq 0$,

$$0 = \lambda_1 t - \mu_1 \bar{T}_1(t), \qquad 0 = \lambda_2 t - \mu_2 \bar{T}_2(t), \qquad 0 = \mu_2 \bar{T}_2(t) - \mu_3 \bar{T}_3(t),$$
$$\bar{I}_1(t) = t - \bar{T}_1(t) - \bar{T}_2(t), \qquad \bar{I}_2(t) = t - \bar{T}_3(t).$$

The result now follows on recalling the definition of $\bar{T}^*(\cdot)$ and Assumption 2.2. $\square$

**Acknowledgments.** We like to thank two referees for careful and in depth reports. Their numerous corrections and suggestions have significantly improved this paper.

Department of Statistics
University of North Carolina
Chapel Hill, North Carolina 27599-3260
USA
e-mail: budhiraj@email.unc.edu
e-mail: apghosh@email.unc.edu